\documentclass[11pt]{amsart}
\usepackage{amsmath, amssymb, amsthm, amsrefs}
\usepackage{latexsym}
\usepackage{indentfirst}
\usepackage{graphicx}
\usepackage{placeins}
\usepackage{booktabs}
\usepackage{algorithm}
\usepackage{algorithmic}
\usepackage{multirow}

\theoremstyle{plain}

\theoremstyle{definition}

\theoremstyle{remark}

\DeclareMathOperator{\tr}{tr}

\DeclareMathOperator{\argmin}{arg min}

\newcommand{\average}[1]{\left\langle#1\right\rangle}

\newcommand{\bbm}{\begin{bmatrix}}
\newcommand{\ebm}{\end{bmatrix}}

\newcommand{\p}{\partial}

\begin{document}

\title[Tensor Network Skeletonization]{Tensor Network Skeletonization}

\author{Lexing Ying} 

\address{
  Department of Mathematics and Institute for Computational and Mathematical Engineering,
  Stanford University,
  Stanford, CA 94305
}

\email{lexing@stanford.edu}

\thanks{This work was partially supported by the National Science
  Foundation under award DMS-1521830 and the U.S. Department of
  Energy’s Advanced Scientific Computing Research program under award
  DE-FC02-13ER26134/DE-SC0009409. The author thanks Lin Lin,
  Xiao-liang Qi, Tao Xiang, Zhiyuan Xie, and Wotao Yin for stimulating
  discussions on various parts of this work.}



\begin{abstract}
  We introduce a new coarse-graining algorithm, tensor network
  skeletonization, for the numerical computation of tensor
  networks. This approach utilizes a structure-preserving
  skeletonization procedure to remove short-range correlations
  effectively at every scale. This approach is first presented in the
  setting of 2D statistical Ising model and is then extended to higher
  dimensional tensor networks and disordered systems. When applied to
  the Euclidean path integral formulation, this approach also gives
  rise to new efficient representations of the ground states for 1D
  and 2D quantum Ising models.
\end{abstract}

\maketitle

%-----------------------------------
\section{Introduction}\label{sec:intro}
This paper is concerned with the numerical computation of tensor
networks (see \cite{Orus2014} for a good introduction for tensor
networks). Recently, tensor networks have received a lot of attention
in computational statistical mechanics and quantum mechanics as
they offer a convenient and effective framework for representing
\begin{itemize}
\item the partition functions for the classical spin systems in
  statistical mechanics and
\item the ground and thermal states of quantum many body systems
  through the Euclidean path integral formulation.
\end{itemize}

\subsection{Definition}
A tensor network is associated with a triple $(V,E,\{T^i\}_{i\in V})$
where $V$ is a set of vertices, $E$ is a set of edges, and $T^i$ is a
tensor at vertex $i\in V$.
\begin{itemize}
\item The degree of the vertex $i\in V$ is denoted as $d_i$.
\item Each edge $e\in E$ is associated with a bond dimension $\chi_e$.
\item For each vertex $i\in V$, $T^i$ is a $d_i$-tensor. Each of the
  $d_i$ components of $T^i$ is associated with one of the adjacent
  edges of the vertex $i$ and the dimension of this component is equal
  to the bond dimension of the associated edge.
\end{itemize}

The edge set $E$ can often be partitioned as the disjoint union $E =
E_I \cup E_B$, where $E_I$ is the set of the interior edges that link
two vertices in $V$ and $E_B$ is the set of the boundary edges with
one endpoint in $V$ and the other one open. Once the triple
$(V,E,\{T^i\}_{i\in V})$ is specified, the tensor network represents a
tensor that is obtained via contracting all interior edges in $E_I$.
The result is an $|E_B|$-tensor and is denoted as
\begin{equation}\label{eq:tr}
  \tr_{E_I} \left( \bigotimes_{i\in V} T^i \right).
\end{equation}
When the set $E_B$ is empty, the tensor network contracts to a
scalar. Throughout this paper, we follow the following notation
conventions.
\begin{itemize}
\item The lower-case letters such as $i$ and $j$ are used to denote
  vertices in $V$.
\item The lower-case letters such as $a$, $b$, $c$, $d$, $e$, and $f$
  are used to denote the edges in $E$.
\item The upper-case letters such as $T$, $U$, and $V$ are used to
  denote tensors.
\end{itemize}

\begin{figure}[h!]
  \includegraphics[scale=0.2]{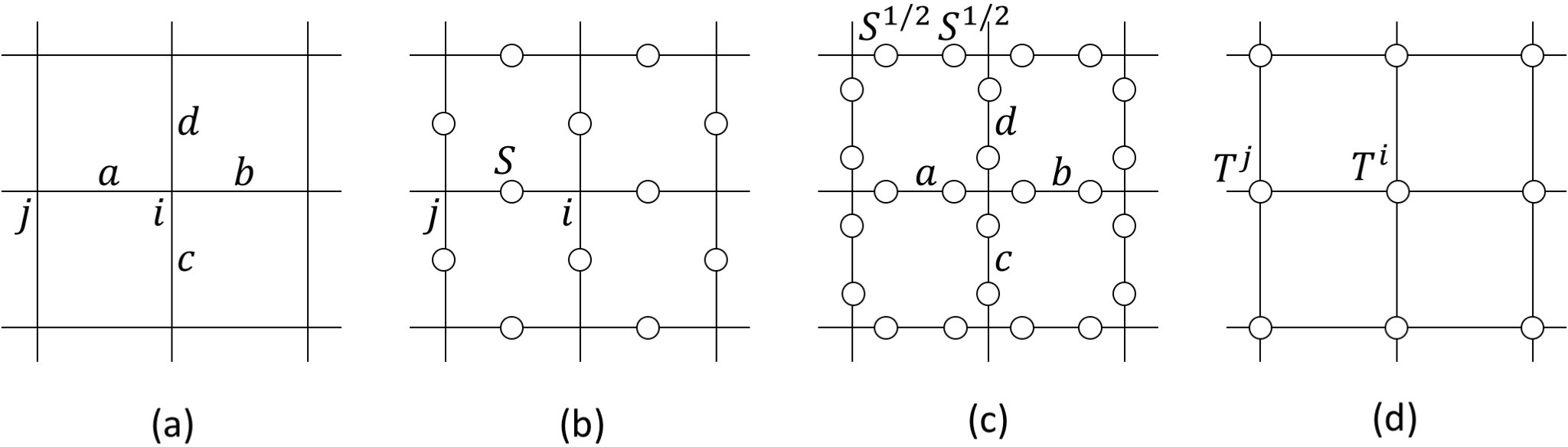}
  \caption{Representing the partition function $Z_N(\beta)$ of 2D
    statistical Ising model using a tensor network. (a) The vertices
    and edges of the tensor network. (b) The 2-tensor $S$ associated
    with each edge. (c) Introducing $S^{1/2}$ splits $S$ into the
    product of two 2-tensors. (d) Contracting the four 2-tensors
    adjacent to a vertex $i$ forms the 4-tensor $T^i$.
  }
  \label{fig:spin2tensor}
\end{figure}

%-------
The framework of tensor networks is a powerful tool since
mathematically it offers efficient representations for high
dimensional functions or probability distributions with certain
underlying geometric structures. For example, let us consider the 2D
statistical Ising model on a periodic square lattice. The vertex set
$V$ consists of the lattice points of an $n\times n$ Cartesian grid,
where $N=n\times n$ is the number of vertices. The edge set $E$
consists of the edges between horizontal and vertical neighbors,
defined using the periodic boundary condition (see Figure
\ref{fig:spin2tensor}(a)). Here $|E|=2N$, $E_I=E$, and $E_B =
\emptyset$.

At temperature $T$, the partition function $Z_N(\beta)$ for the
inverse temperature $\beta=1/T$ is given by
\[
Z_N(\beta) = \sum_{\sigma} e^{-\beta H_N(\sigma)}, \quad H_N(\sigma) =
-\sum_{(ij)\in E} \sigma_i \sigma_j,
\]
where $\sigma=(\sigma_1,\ldots,\sigma_N)$ stands for a spin
configuration at $N$ vertices with $\sigma_i=\pm 1$ and the sum of
$\sigma$ is taken over all $2^N$ configurations. Here $(ij)$ is the
edge between two adjacent vertices $i$ and $j$ and the sum in
$H_N(\sigma)$ is over these $2N$ edges.

In order to write $Z_N(\beta)$ in the form of a tensor network, one
approach is to introduce a $2\times 2$ matrix $S$
\[
S = \begin{pmatrix} e^\beta & e^{-\beta} \\ e^{-\beta} & e^{\beta} \end{pmatrix},
\]
which is the multiplicative term in $Z_N(\beta)$ associated with an
edge between any two adjacent vertices. The partition function
$Z_N(\beta)$ is built from the $S$ matrices over all edges in $E$ (see
Figure \ref{fig:spin2tensor}(b)). Since $S$ is a symmetric matrix, its
symmetric square root $S^{1/2}$ is well defined with the following
element-wise identity:
\[
S_{ij} = \sum_a S^{1/2}_{ia} S^{1/2}_{aj}
\]
where $a$ denotes the edge that connects $i$ and $j$ (see Figure
\ref{fig:spin2tensor}(c)). Here and throughout the paper, the
lower-case letters (e.g. $i$, $j$, and $k$) for denoting a vertex in
$V$ are also used for the running index associated with that
vertex. The same applies to the edges: the lower-case letters
(e.g. $e$ and $f$) for denoting an edge in $E$ are also used as the
running index associated with that edge.

At each vertex $i$, one can then introduce a $4$-tensor $T^i$
\begin{equation}\label{eq:Ti}
T^i_{abcd} = \sum_i S^{1/2}_{ia}S^{1/2}_{ib}S^{1/2}_{ic}S^{1/2}_{id},
\end{equation}
which essentially contracts the four $S^{1/2}$ tensors adjacent to the
vertex $i$ (see Figure \ref{fig:spin2tensor}(c)). Finally, the
partition function $Z_N(\beta)$ can be written as
\[
Z_N(\beta) = \tr_{E} \left( \bigotimes_{i\in V} T^i \right)
\]
(see Figure \ref{fig:spin2tensor}(d)).

%-------
\subsection{Previous work}

One of the main computational tasks is how to evaluate tensor networks
accurately and efficiently. The naive tensor contraction following the
definition \eqref{eq:tr} is computationally prohibitive since its
running time grows exponentially with the number of vertices.

In recent years, there has been a lot of work devoted to efficient
algorithms for evaluating tensor networks. In \cite{Levin2007}, Levin
and Nave introduced the tensor renormalization group (TRG) as probably
the first practical algorithm for this task. When applied to the 2D
statistical Ising models, this method utilizes an alternating sequence
of tensor contractions and singular value decompositions. However, one
problem with TRG is the accumulation of short-range correlation, which
increases the bond dimensions and computational costs dramatically as
the method proceeds.

In a series of papers \cite{Xie2009,Xie2012,Zhao2010}, Xiang et al
introduced the higher order tensor renormalization group (HOTRG) as an
extension of TRG to address 3D classical spin systems. The same group
has also introduced the second renormalization group (SRG) based on
the idea of approximating the environment of a local tensor before
performing reductions. SRG typically gives more accurate
results. However, the computation time of SRG tends to grow
significantly with the size of the local environment and it is also
not clear how to generalize this technique to systems that are not
translationally invariant.

In \cite{Gu2009}, Gu and Wen introduced the method of tensor
entanglement filtering renormalization (TEFR) as an improvement of TRG
for 2D systems. Comparing with TRG, this method makes an extra effort
in removing short-range correlations and hence produces more accurate
and efficient computations.

More recently in \cite{Evenbly2015,Evenbly2015B}, Evenbly and Vidal
proposed the tensor network renormalization (TNR). The key step of TNR
is to apply the {\em disentanglers} to remove short-range
correlation. These disentanglers appeared earlier in the work of the
multiscale entanglement renormalization ansatz (MERA)
\cite{Vidal2008}. For a fixed bond dimension, TNR gives significantly
more accurate results compared to TRG, but at the cost of increasing
the computational complexity. However, it is not clear how to extend
the approach of TNR to systems in higher dimensions.

These approaches have significantly improved the efficiency and
accuracy of the computation of tensor networks. From a computational
point of view, it would be great to have a general algorithm that have
the following three properties:
\begin{itemize}
\item removing the short-range correlation efficiently in order to
  keep bond dimension and computational cost under control, and
\item extending to 3D and 4D tensor networks, and 
\item extending to systems that are not translationally invariant,
  such as disordered systems.
\end{itemize}
However, as far as we know, none of these methods achieves all three
properties simultaneously.

%-------
\subsection{Contribution and outline}

Building on top of the previous work in the physics literature, we
introduce a new coarse-graining approach, called the tensor network
skeletonization (TNS), as a first step towards building such a general
algorithm. At the heart of this approach is a new procedure called the
{\em structure-preserving skeletonization}, which removes short-range
correlation efficiently while maintaining the structure of a local
tensor network. This allows us to generalize TNS quite
straightforwardly to spin systems of higher dimensions. In addition,
we also provide a simple and efficient algorithm for performing the
structure-preserving skeletonization. This allows for applying TNS
efficiently to systems that are not translationally invariant.

The rest of this paper is organized as follows. Section
\ref{sec:intro} summarizes the basic tools used by the usual tensor
network algorithms and introduces the structure-preserving
skeletonization. Section \ref{sec:2D} is the main part of the paper
and explains TNS for 2D statistical Ising model. Section \ref{sec:3D}
extends the algorithm to 3D statistical Ising model. Section
\ref{sec:ground} discusses how to build efficient representations of
the ground states of 1D and 2D quantum Ising models using
TNS. Finally, Section \ref{sec:conc} discusses some future work.

%-----------------------------------
\section{Basic tools}\label{sec:basic}

\subsection{Local replacement}\label{sec:basiclr}
\begin{figure}[h!]
  \includegraphics[scale=0.2]{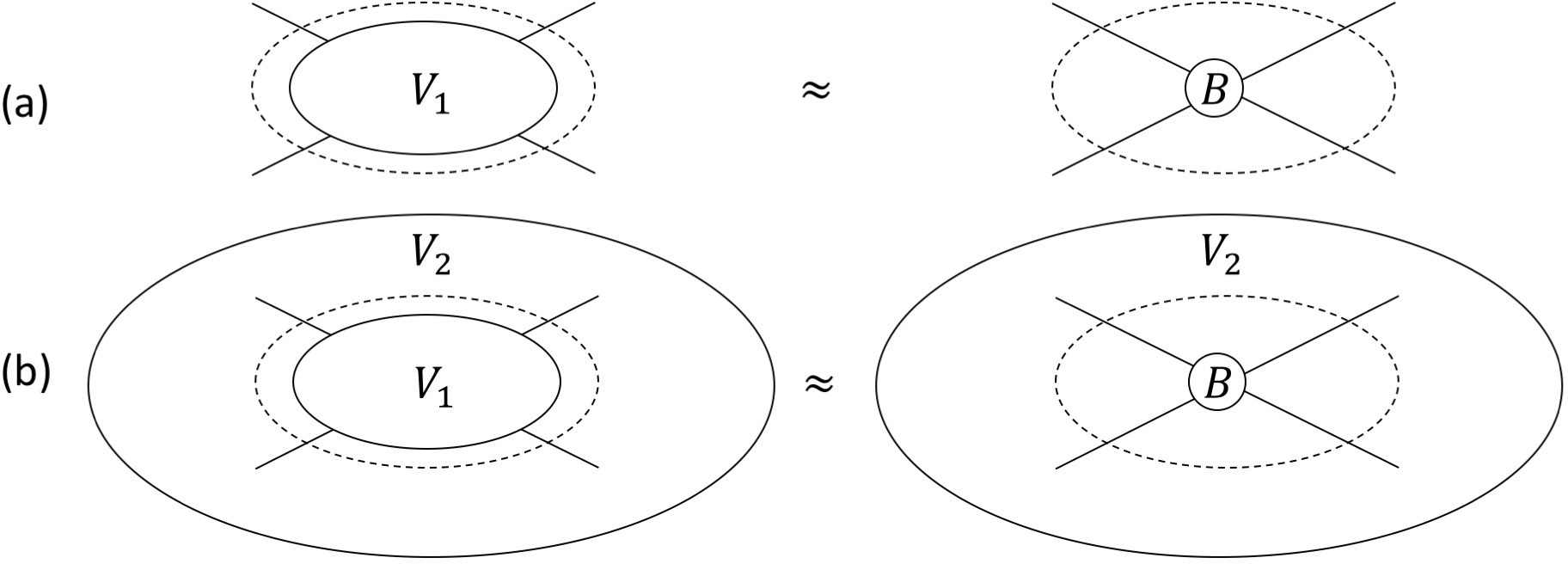}
  \caption{Local replacement. (a) Part of the tensor network
    associated with vertices in $V_1$ is approximated by a simplified
    tensor network $B$.  (b) Locally replacing $V_1$ with $B$ results
    in a approximation of the whole tensor network.}
  \label{fig:Bapprox}
\end{figure}

The basic building blocks of all tensor network algorithms are {\em
  local replacements}. Suppose that vertex $V$ and edge set $E$ of a
tensor network $(V,E,\{T^i\}_{i\in V})$ are partitioned as follows
\[
V = V_1 \cup V_2, \quad E = E_1 \cup E_2 \cup E_{12},
\]
where $E_1$ and $E_2$ are the sets of interior edges of $V_1$ and
$V_2$, respectively, and $E_{12}$ is the set of edges that link across
$V_1$ and $V_2$. Such a partition immediately gives an identity
\begin{equation}\label{eq:lclrepid}
  \tr_{E} \left(\bigotimes_{i\in V}  T^i\right) = \tr_{E_2\cup E_{12}}
  \left(  \left(\bigotimes_{i\in V_2}T^i\right) \bigotimes
  \tr_{E_1}\left(\bigotimes_{i\in V_1}T^i\right) \right).
\end{equation}
Assume now that there exists another tensor network $B$ for which the
following approximation holds
\[
B \approx \tr_{E_1}\left(\bigotimes_{i\in V_1} T^i\right)
\]
(see Figure \ref{fig:Bapprox}(a)). Typically $B$ is much simpler in
terms of the number of the vertices and/or the bond dimensions of the
edges. A {\em local replacement} refers to replacing
$\tr_{E_1}\left(\bigotimes_{i\in V_1} T^i\right)$ in
\eqref{eq:lclrepid} with $B$ to get a simplified approximation
\[
\tr_{E} \left(\bigotimes_{i\in V}  T^i\right)
\approx \tr_{E_2\cup E_{12}} \left( \left(\bigotimes_{i\in V_2}T^i\right) \bigotimes B \right)
\]
(see Figure \ref{fig:Bapprox}(b)). Most algorithms for tensor networks
apply different types of local replacements successively until the
tensor network is simplified to a scalar or left with only boundary
edges.

\begin{figure}[h!]
  \includegraphics[scale=0.2]{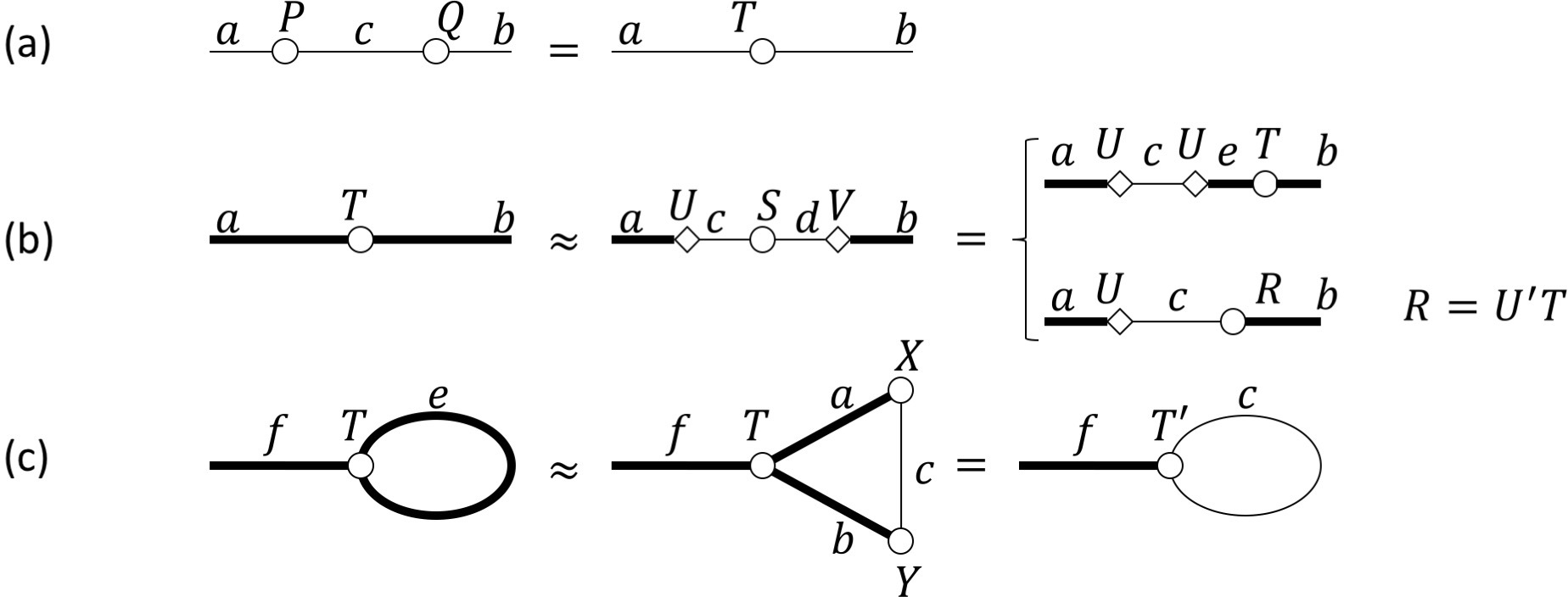}
  \caption{Instances of local replacements. (a) Contraction. (b)
    Projection. (c) Structure-preserving skeletonization.  }
  \label{fig:Localtools}
\end{figure}

The simplest instance of local replacement is the {\em tensor
  contraction} and it simply combines two adjacent tensors into a
single one. For example, let $P$ be a 2-tensor adjacent to edges $a$
and $c$ and $Q$ be another 2-tensor adjacent to edges $c$ and $b$ (see
Figure \ref{fig:Localtools}(a)). The resulting 2-tensor $T$ obtained
from contracting $P$ and $Q$ is simply the product of $P$ and $Q$,
i.e.,
\[
T_{ab} = \sum_{c}P_{ac}Q_{cb},
\]
(see Figure \ref{fig:Localtools}(a)). Often when the contraction is
applied, the edges $a$ and $b$ typically come from grouping a set of
multiple edges.

A second instance is called the {\em projection}. Typically it is
carried out by performing a singular value decomposition of $T$
followed by thresholding small singular value, i.e.,
\[
T \approx U S V', \quad T_{ab} \approx \sum_{cd} U_{ac} S_{cd} V_{bd}
\]
(see Figure \ref{fig:Localtools}(b)). Here $U$ and $V$ are both
orthogonal matrices and $S$ is a diagonal matrix.  Due to the
truncation of small singular values, the bond dimensions at edges $c$
and $d$ can be significantly smaller compared to the ones of $a$ and
$b$.  Throughout this paper, each orthogonal matrix shall be denoted
by a diamond in the figures.  As with the contraction, each of the
indices $a$ and $b$ often comes from grouping a set of multiple edges.
The SVD-based projection can also be modified slightly to a few
equivalent forms (see Figure \ref{fig:Localtools}(b))
\begin{align*}
  & T \approx UU'T, \quad T_{ab}  \approx \sum_{ce} U_{ac}U_{ec}T_{eb}\\
  & T \approx UR,   \quad R=U'T, \quad T_{ab} \approx \sum_{c} U_{ac} R_{cb}.
\end{align*}
In the rest of this paper, we refer to the first one as the
$UU'T$-projection and the second one as the $UR$-projection.

Another instance of local replacements uses the {\em disentanglers}
introduced in \cite{Vidal2008} and it plays a key role in the work of
tensor network renormalization (TNR) \cite{Evenbly2015} as mentioned
above. Since the tensor network skeletonization (TNS) approach of this
paper does not depend on the disentanglers, we refer to the references
\cite{Evenbly2015,Vidal2007,Evenbly2009} for more detailed discussions
of them.

\subsection{Structure-preserving skeletonization}\label{sec:basicsps}

At the center of the TNS approach is a new type of local replacement
called the {\em structure-preserving skeletonization}. The overall
goal is to reduce the bond dimensions of the interior edges of a loopy
local tensor network without changing its topology. In the simplest
setting, consider a 3-tensor $T$ with two of its components marked
with a same edge $e$ (and thus to be contracted). The
structure-preserving skeletonization seeks two 2-tensors $X$ and $Y$
(see Figure \ref{fig:Localtools}(c)) such that
\begin{equation}\label{eq:sps}
  \tr_e T \approx \tr_{abc} (X\otimes T \otimes Y),\quad
  \sum_e T_{eef} \approx \sum_{abc} X_{ac} T_{abf} Y_{bc},
\end{equation}
and also the bond dimension of edge $c$ should be significantly
smaller compared to the bond dimension of edge $e$. This is possible
because there might exists short-range correlations within the loop
that can be removed from the viewpoint of the exterior of this local
tensor network.

A convenient way to reformulate the problem is to view $T$ as a
$\chi_e\times\chi_e$ matrix with each entry $T_{ab}$ equal to a
$\chi_f$-dimensional vector and view $X$ and $Y$ as matrices. Then one
can rewrite the condition in \eqref{eq:sps} as
\begin{equation}\label{eq:spsmat}
\tr_e T \approx \tr_c(X^* T Y)
\end{equation}
where the products between $X$, $T$, and $Y$ are understood as matrix
multiplications.

As far as we know, there does not exist a simple and robust numerical
linear algebra routine that solve this approximation problem
directly. Instead, we propose to solve the following regularized
optimization problem
\[
\min_{X,Y} \|\tr_e T - \tr_c (X^* T Y) \|_2^2 + \alpha
(\|X\|_F^2+\|Y\|_F^2),
\]
where the constant $\alpha$ is a regularization parameter and is
typically chosen to be sufficiently small. This optimization problem
is non-convex, however it can be solved effectively in practice using
the alternating least square algorithm once a good initial guess is
available. More precisely, given a initial guess for $X^{(0)}$ and
$Y^{(0)}$, one alternates the following two steps for $n=0,1,\ldots$
until convergence
\begin{align*}
  X^{(n+1)} &= \argmin_{X} \|\tr_e T - \tr_{c} (X^* T  Y^{(n)}) \|_2^2 + \alpha \|X\|_F^2\\
  Y^{(n+1)} &= \argmin_{Y} \|\tr_e T - \tr_{c} ((X^{(n+1)})^* T  Y) \|_2^2 + \alpha \|Y\|_F^2.
\end{align*}
Since each of the two steps is a least square problem in $X$ or $Y$,
they can be solved efficiently with standard numerical linear algebra
routines. The numerical experience shows that, when starting from
well-chosen initial guesses, this alternating least square algorithm
converges after a small number of iterations to near optimal
solutions.

%-----------------------------------
\section{TNS for 2D statistical Ising models}\label{sec:2D}

We start with a 2D statistical Ising model on an $n\times n$ lattice
with the periodic boundary condition. Following the discussion in
Section \ref{sec:intro}, we set the vertex set $V_0$ to be an $n\times
n$ Cartesian grid. Each vertex $i$ is identified with a tuple
$(i_1,i_2)$ with $i_1,i_2\in [n]=\{0,1,\ldots,n-1\}$. The edge set
$E_0$ consists of the edges between horizontal and vertical neighbors
of the Cartesian grid modulus periodicity. This setup also gives rise
to an $n\times n$ array of plaquettes. If a plaquette has vertex
$i=(i_1,i_2)$ at its lower-left corner, then we shall index this
plaquette with $i=(i_1,i_2)$ as well. Here $N=n^2$ is the total number
of spins and we assume without loss of generality that $n=2^L$.

%-----------
\subsection{Partition function}\label{sec:2Dpf}

Following the discussion in Section \ref{sec:intro}, the partition
function $Z_N(\beta)$ can be represented using a tensor network
$(V^0,E^0,\{T^i\}_{i\in V_0})$ where $T^i$ are given in \eqref{eq:Ti}.
Let $\chi$ be a predetermined upper bound for the bond dimension of
the edges of the tensor network. One can assume without loss of
generality that the bond dimension $\chi_e$ for the edge $e\in E^0$ is
close to this constant $\chi$. When $\chi$ is significantly larger
than 2, this can be achieved by contracting each $2\time 2$
neighborhood of tensors into a single tensor. For example, when
$\chi=4$, one round of such contractions brings $\chi_e=\chi=4$.

%----
\subsubsection{Algorithm}\label{sec:2Dpfal}

\begin{figure}[h!]
  \includegraphics[scale=0.2]{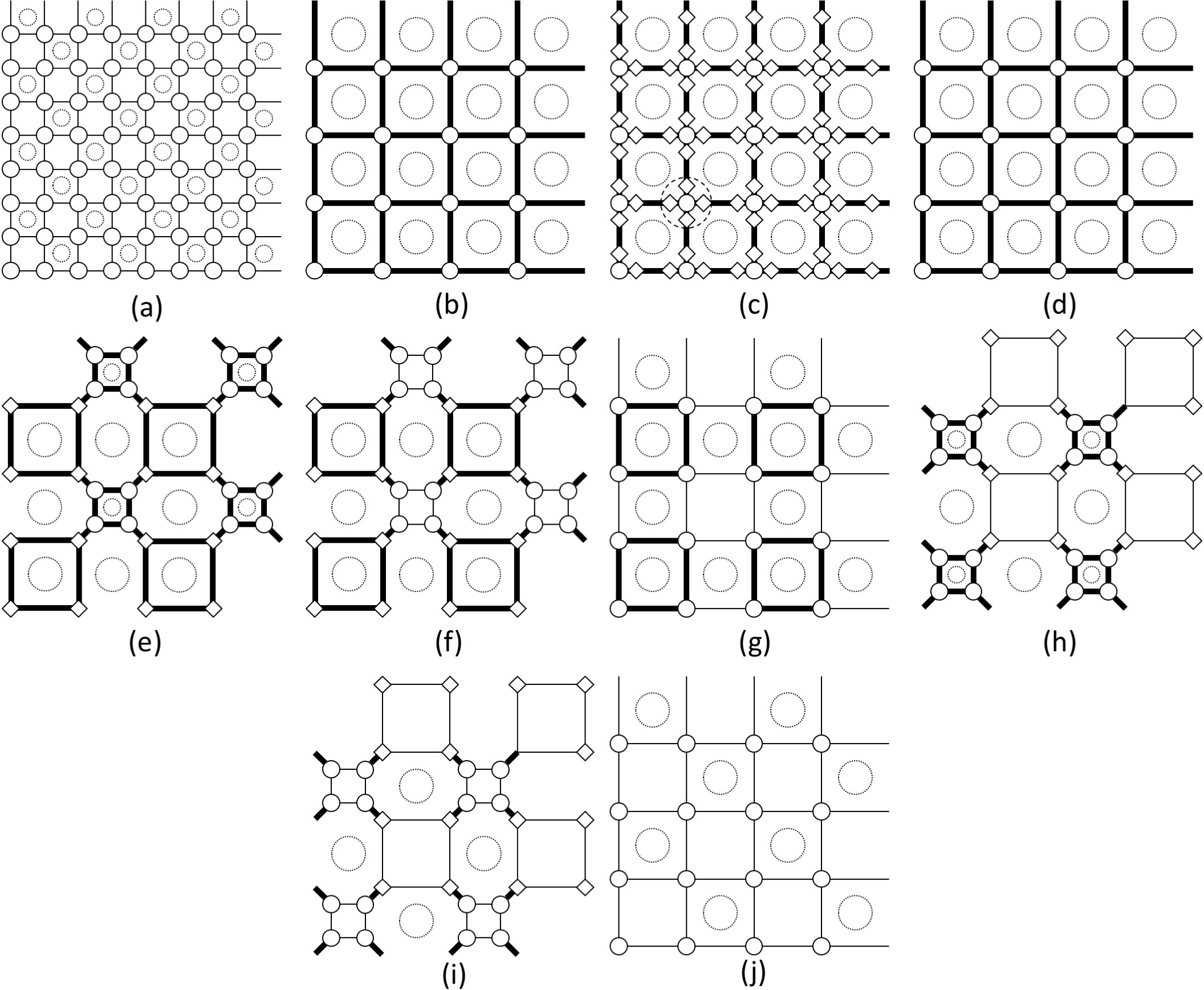}
  \caption{A single iteration of the tensor network skeletonization
    algorithm. The starting point is a tensor network with bond
    dimension $\chi$ and short-range correlation removed in the
    $(0,0)_2$ and $(1,1)_2$ plaquettes. The final point is a
    coarse-grained tensor-network with 1/4 vertices (or tensors). This
    coarse-grained tensor also has bond dimensions equal to $\chi$ and
    has short-range correlation removed for the (larger) $(0,0)_2$ and
    $(1,1)_2$ plaquettes. The bold lines stand for edges with bond
    dimensions equal to $O(\chi^2)$. }
  \label{fig:2D}
\end{figure}

The TNS algorithm consists of a sequence of coarse-graining
iterations. At the beginning of the $\ell$-th iteration, one holds a
tensor network $(V^\ell,E^\ell,\{T^i\}_{i\in V_\ell})$ at level $\ell$
with $(n/2^\ell)\times (n/2^\ell)$ vertices. With the exception of the
$0$-th iteration, we require that the following iteration invariance
to hold:
\begin{itemize}
\item for each plaquette with index equal to $(0,0)$ or $(1,1)$
  modulus 2, the short-range correlation within this plaquette has
  already been eliminated (see Figure \ref{fig:2D}(a)).
\end{itemize}
In what follows, we refer to those plaquettes with index equal to
$(0,0)$ modulus 2 as $(0,0)_2$ plaquettes and similarly those with
index equal to $(1,1)$ modulus 2 as $(1,1)_2$ plaquettes. In Figure
\ref{fig:2D}(a), the dotted circles denote the existence of
short-range correlation. Notice that these circles do not appear in
the $(0,0)_2$ plaquettes and $(1,1)_2$ plaquettes.  The $\ell$-th
iteration consists of the following steps.
\begin{enumerate}
\item Merge the tensors at the four vertices of each $(0,0)_2$
  plaquette into a single tensor (see Figure \ref{fig:2D}(b)). This
  requires a couple of contractions defined in Section
  \ref{sec:basic}. The $(1,1)_2$ plaquettes are stretched and this
  results a new graph that contains only $1/4$ of vertices. The
  tensors at the new vertices are identical but the bond dimension of
  the new edges are equal to $\chi^2$ (shown using the bold
  lines). Since these $(1,1)_2$ plaquettes at level $\ell$ do not
  contain short-range correlation at level $\ell$, no short-range
  correlation of level $\ell$ will survive at level $\ell+1$. However,
  there are new short-range correlation of level $\ell+1$ in the
  tensor network and these are marked with larger dotted circles
  inside the new plaquettes at level $\ell+1$. The key task for the
  rest of the iteration is to remove part of these short-range
  correlations at level $\ell+1$ and reduce the bond dimension from
  $\chi^2$ back to $\chi$ at the same time.
\item For each vertex $i$ in Figure \ref{fig:2D}(b), denote the tensor
  at $i$ by $T^i_{abcd}$ where $a$, $b$, $c$, and $d$ refer to the
  left, right, bottom, and top edges. Applying two $UU'T$-projections
  (the first one with the left edge vs. the rest and the second one
  with for the bottom edge vs. the rest) to $T^i_{abcd}$ effectively
  inserts two orthogonal (diamond) matrices in each of these edges
  (see Figure \ref{fig:2D}(c)). Most TRG algorithms utilize this step
  to reduce the bond dimension directly from $\chi^2$ back to $\chi$
  and thus incurring a significant error. In TNS however, the bond
  dimension at this point is kept close to $\chi^2$. This introduces a
  much smaller truncation error as compared to the TRG algorithms.
\item At each vertex $i$ in Figure \ref{fig:2D}(c), merge the tensor
  $T^i_{abcd}$ with the four adjacent orthogonal (diamond) matrices
  (see Figure \ref{fig:2D}(d)). Though the tensor network obtained
  after this step has the same topology as the one in Figure
  \ref{fig:2D}(b), the bond dimension is somewhat reduced and one
  prepares the tensor network for the structure-preserving
  skeletonization.
\item For each $(1,1)_2$ plaquette in Figure \ref{fig:2D}(d), apply
  the $UR$-projection to the $4$-tensor at each of its corners. Here
  the two edges adjacent to the plaquette are grouped together. Notice
  that the (round) $R$ tensors are placed close to the $(1,1)_2$
  plaquette while the (diamond) $U$ tensors are placed away from
  it. Though this projection step does not reduce bond dimensions, it
  allows us to treat each $(1,1)_2$ plaquette separately. The
  resulting graph is given in Figure \ref{fig:2D}(e).
\item In this key step, apply the {\em structure-preserving
  skeletonization} to each $(1,1)_2$ plaquette. The details of this
  procedure will be provided below in Section \ref{sec:2Dpfsps}. The
  resulting $(1,1)_2$ plaquette has its short-range correlation
  removed and the bond dimensions of its four surrounding edges are
  reduced from $\chi^2$ to $\chi$ (see Figure \ref{fig:2D}(f)).
\item For each $(1,1)_2$ plaquette, contract back the $UR$-projections
  at each of its four corners. Notice that, due to the
  structure-preserving skeletonizations, the new $R$ tensors have bond
  dimensions equal to $\chi$. The resulting tensor network (see Figure
  \ref{fig:2D}(g)) is similar to the one in Figure \ref{fig:2D}(d) but
  now the short-range correlations in the $(1,1)_2$ plaquettes are all
  removed.
\item Now repeat the previous three steps to the $(0,0)_2$
  plaquettes. This is illustrated in Figure \ref{fig:2D}(h), (i), and
  (j). The resulting tensor network now has short-range correlation
  removed in both $(0,0)_2$ and $(1,1)_2$ plaquettes. In addition, the
  bond dimension of the edges is reduced back to $\chi$ from $\chi^2$.
\end{enumerate}
This finishes the $\ell$-th iteration. At this point, one obtains a
new tensor network denoted by $(V^{\ell+1},E^{\ell+1},\{T^i\}_{i\in
  V_{\ell+1}})$ that is a self-similar and coarse-grained version of
$(V^\ell,E^\ell,\{T^i\}_{i\in V_\ell})$. Since the short-range
correlations in both $(0,0)_2$ and $(1,1)_2$ plaquettes are removed,
this new tensor network satisfies the iteration invariance and it can
serve as the starting point of the $(\ell+1)$-th iteration.

Following this process, the TNS algorithm constructs a sequence of
tensor networks
\[
(V^\ell,E^\ell,\{T^i\}_{i\in V_\ell}), \quad \ell=0,1,2,\ldots,L.
\]
The last one is a single $4$-tensor with the left and right edges
identified and similarly with the bottom and top edges
identified. Contracting this final tensor gives a scalar value for the
partition function $Z_N(\beta)$.

%----
\subsubsection{Structure-preserving skeletonization}\label{sec:2Dpfsps}

In the description of the algorithm in Section \ref{sec:2Dpfal}, the
missing piece is how to carry out the structure-preserving
skeletonization in order to remove the short-range correlation of a
$(1,1)_2$ plaquette and reduce the bond dimension of its four
surrounding edges (from Figure \ref{fig:2D}(e) to Figure
\ref{fig:2D}(f)).

\begin{figure}[h!]
  \includegraphics[scale=0.2]{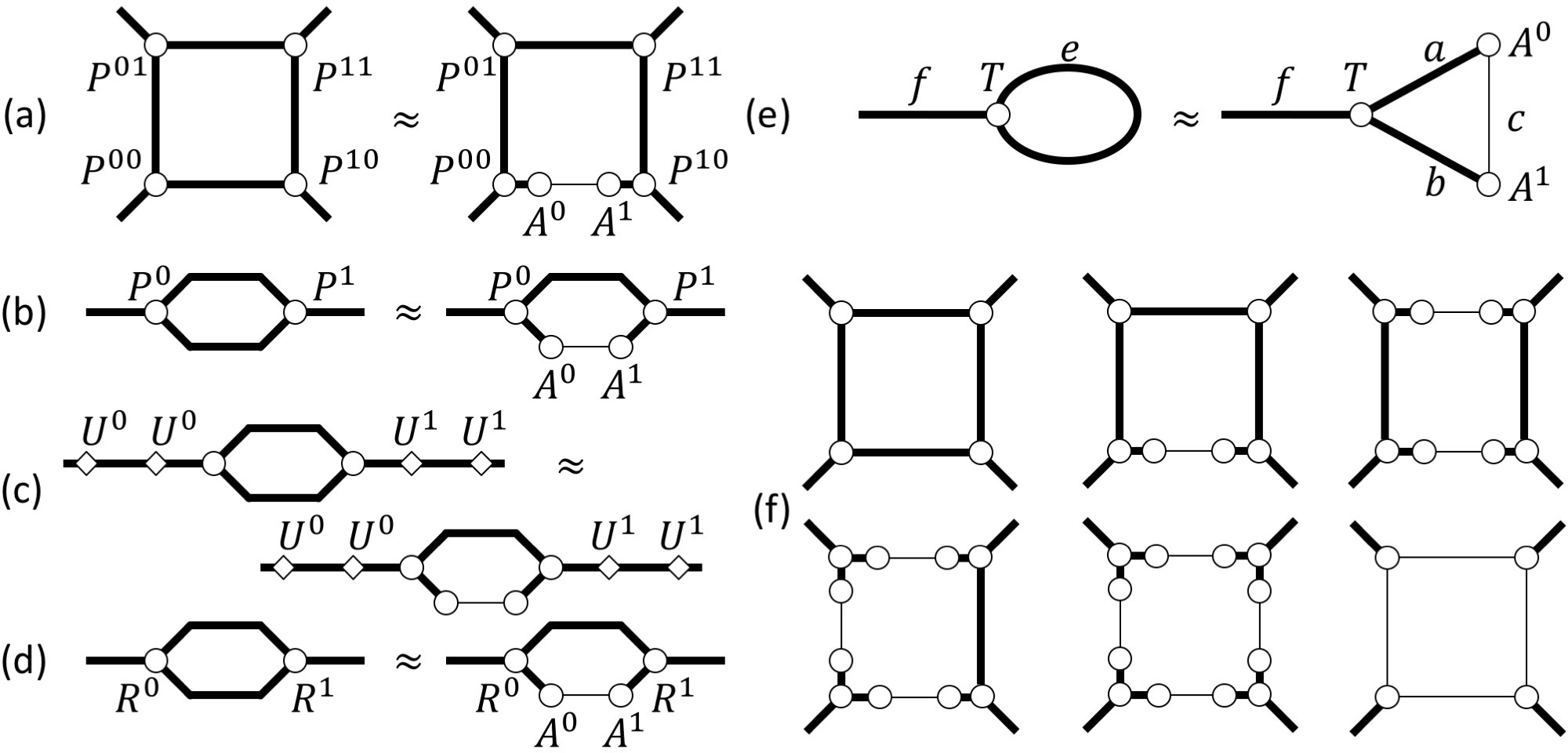}
  \caption{ The structure-preserving skeletonization procedure removes
    the short-range correlation within a $(1,1)_2$ (or $(0,0)_2$)
    plaquette. The bold lines stand for edges with bond dimensions
    $\ge \chi^2$.  }
  \label{fig:2Dskl}
\end{figure}

This procedure is illustrated in Figure \ref{fig:2Dskl} with the four
corner tensors denoted by $P^{00}$, $P^{10}$, $P^{01}$, and
$P^{11}$. Instead of replacing the four corner 3-tensors
simultaneously, this procedure considers the 4 interior edges one by
one and insert for each edge two tensors of size $\chi^2 \times \chi$.

\begin{enumerate}
\item Starting from the bottom edge, we seek two 2-tensors $A^0$ and
  $A^1$ of size $\chi^2 \times \chi$ under the condition that the
  4-tensor represented by the new $(1,1)_2$-plaquette (after inserting
  $A^0$ and $A^1$) approximates the 4-tensor represented by the original
  plaquette (see Figure \ref{fig:2Dskl}(a)).
\item Merge the two left tensors $P^{00}$ and $P^{01}$ into a 3-tensor
  $P^0$ and merge the two right tensors $P^{10}$ and $P^{11}$ into a
  3-tensor $P^1$. After that, the condition is equivalent to the one
  given in Figure \ref{fig:2Dskl}(b). Notice that the two boundary
  edges have bond dimension equal to $\chi^4$.
\item Since the bond dimensions of the two edges between $P^0$ and
  $P^1$ will eventually be reduced to $\chi$, this implies that the
  bond dimensions of the two boundary edges can be cut down to
  $\chi^2=\chi\times\chi$ without affecting the accuracy. For this, we
  perform the $UU'T$-projection to both $P^0$ and $P^1$. This gives
  rise the condition in Figure \ref{fig:2Dskl}(c).
\item Remove the two tensors $U^0$ and $U^1$ at the two endpoints.
  Merge $U^0$ with $P^0$ to obtain a 3-tensor $R^0$ and similarly
  merge $U^1$ with $P^1$ to obtain a 3-tensor $R^1$. The approximation
  condition can now be written in terms of $R^0$ and $R^1$ as in
  Figure \ref{fig:2Dskl}(d).
\item Finally, contracting the top edge between $R^0$ and $R^1$
  results in a new 3-tensor $T$. The approximation condition now takes
  the form given in Figure \ref{fig:2Dskl}(e). This is now exactly the
  setting visited in Section \ref{sec:basicsps} and can be solved
  efficiently using the alternating least square algorithm proposed
  there.
\end{enumerate}
At this point, two tensors $A^0$ and $A^1$ are successfully inserted
into the bottom edge. One can repeat this process now for the top,
left, and right edges in a sequential order. Once this is done,
merging each of the corner tensors with its two adjacent inserted
tensors gives the desired approximation (see Figure \ref{fig:2Dskl}(f)
for this whole process).

A task of reducing the bond dimensions of the four surrounding edges
of a plaquette has appeared before in the work of tensor entanglement
filtering renormalization (TEFR) \cite{Gu2009}. However, the
algorithms proposed there are different the one described here and the
resulting plaquette was used in a TRG setting that does not extend
naturally to higher dimensions tensor network problems.

%----
\subsubsection{A modified version}\label{sec:2Dpfmod}

In terms of coarse-graining the tensor network, each iteration of the
algorithm in Section \ref{sec:2Dpfal} achieves the important goal of
constructing a self-similar version while keeping the bond dimension
constant (equal to $\chi$) (see Figure \ref{fig:2D}(a) and (j) for
comparison).

\begin{figure}[h!]
  \includegraphics[scale=0.2]{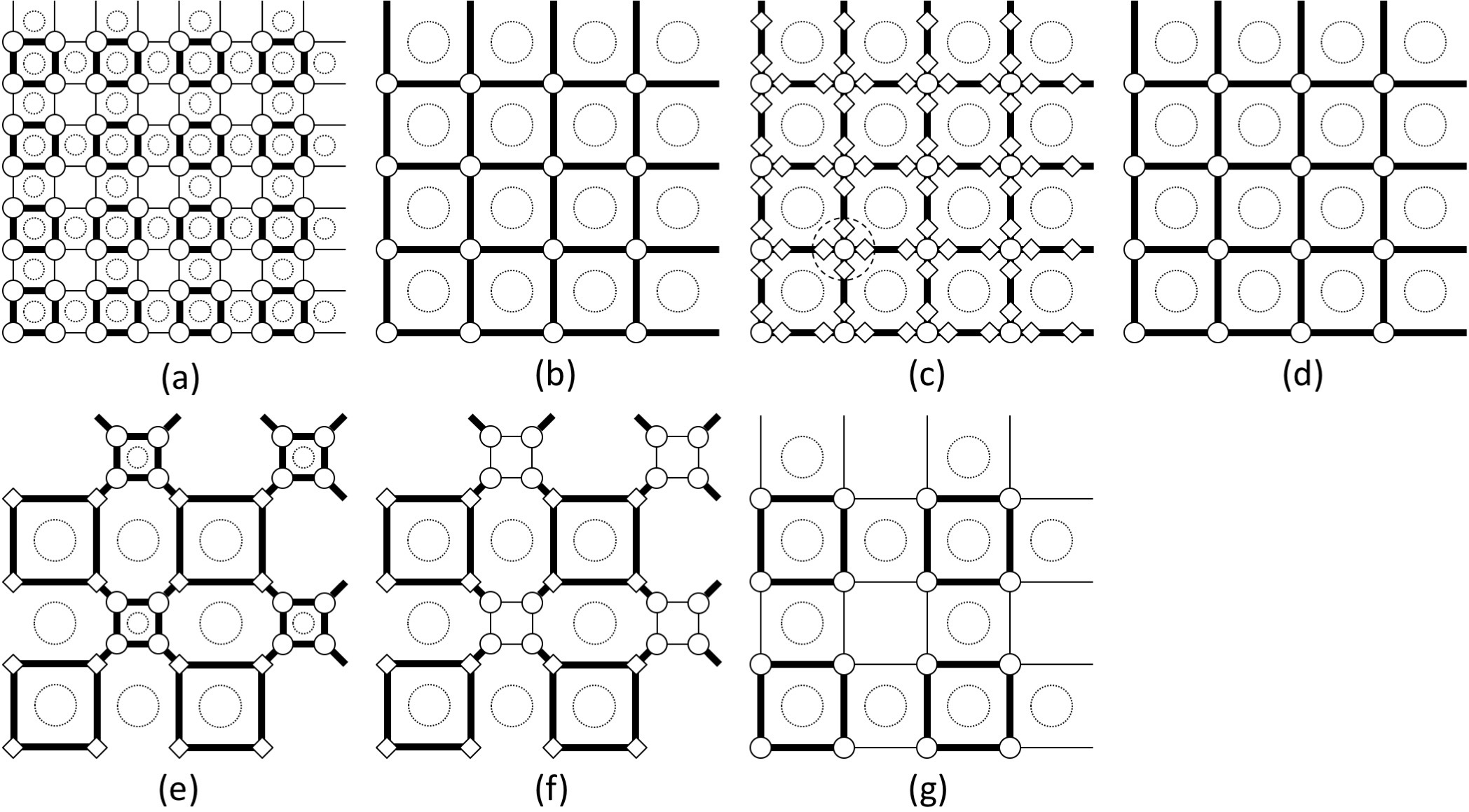}
  \caption{ A single iteration of the modified tensor network
    skeletonization algorithm. The starting point is a tensor network
    with short-range correlation removed in $(1,1)_2$ plaquettes. The
    final point is a coarse-grained tensor-network with 1/4 vertices
    (tensors). This coarse-grained tensor also has bond dimensions
    equal to $\chi$ around the $(1,1)_2$ plaquettes and $\chi^2$
    around the $(0,0)_2$ plaquettes. The short-range correlation is
    removed for the (larger) $(1,1)_2$ plaquettes. The bold lines
    stand for edges with bond dimensions $\chi^2$.  }
  \label{fig:2Dp}
\end{figure}

However, for the purpose of merely computing the partition function
$Z_N(\beta)$, a part of the work is redundant. More specifically, at
the end of the $\ell$-th iteration, the structure-preserving
skeletonization is also performed to the $(0,0)_2$-plaquettes at level
$\ell+1$ to remove their short-range correlations. However, right at
the beginning of the next iteration, a merging step contracts the four
corner tensors of each $(0,0)_2$-plaquette. By eliminating this
structure-preserving skeletonization for the $(0,0)_2$ plaquettes, one
obtains a modified version of the algorithm (see Figure \ref{fig:2Dp})
that can potentially be computationally more efficient.

Compared with the algorithm illustrated in Figure \ref{fig:2D}, the
main differences are listed as follows
\begin{itemize}
\item The iteration invariance is that, at the beginning of each
  iteration, only the short-range correlations of the $(1,1)_2$
  plaquettes are removed.  Therefore, the bond dimensions of the edges
  around a $(0,0)_2$ plaquette are equal to $\chi^2$. This is not so
  appealing from the viewpoint of approximating a tensor network with
  minimal bond dimension. However, as one can see from Figure
  \ref{fig:2Dp}(a) to Figure \ref{fig:2Dp}(b), a contraction step is
  applied immediately to these $(0,0)_2$ plaquettes so that the high
  bond dimensions do not affect subsequent computations.
\item From Figure \ref{fig:2Dp}(e) to Figure \ref{fig:2Dp}(f), the
  structure-preserving skeletonization is only applied to the
  $(1,1)_2$ plaquettes.
\item In Figure \ref{fig:2Dp}(g), the resulting tensor network at
  level $\ell+1$ satisfies the new iteration invariance and hence it
  can serve as the starting point of the next iteration.
\end{itemize}
As we shall see in Section \ref{sec:2Dobal}, this modified algorithm
also has the benefit of incurring minimum modification when evaluating
observables using the impurity method.

%----
\subsubsection{Numerical results}\label{sec:2Dpfna}

\begin{figure}[h!]
  \begin{tabular}{cc}
    \includegraphics[height=2in]{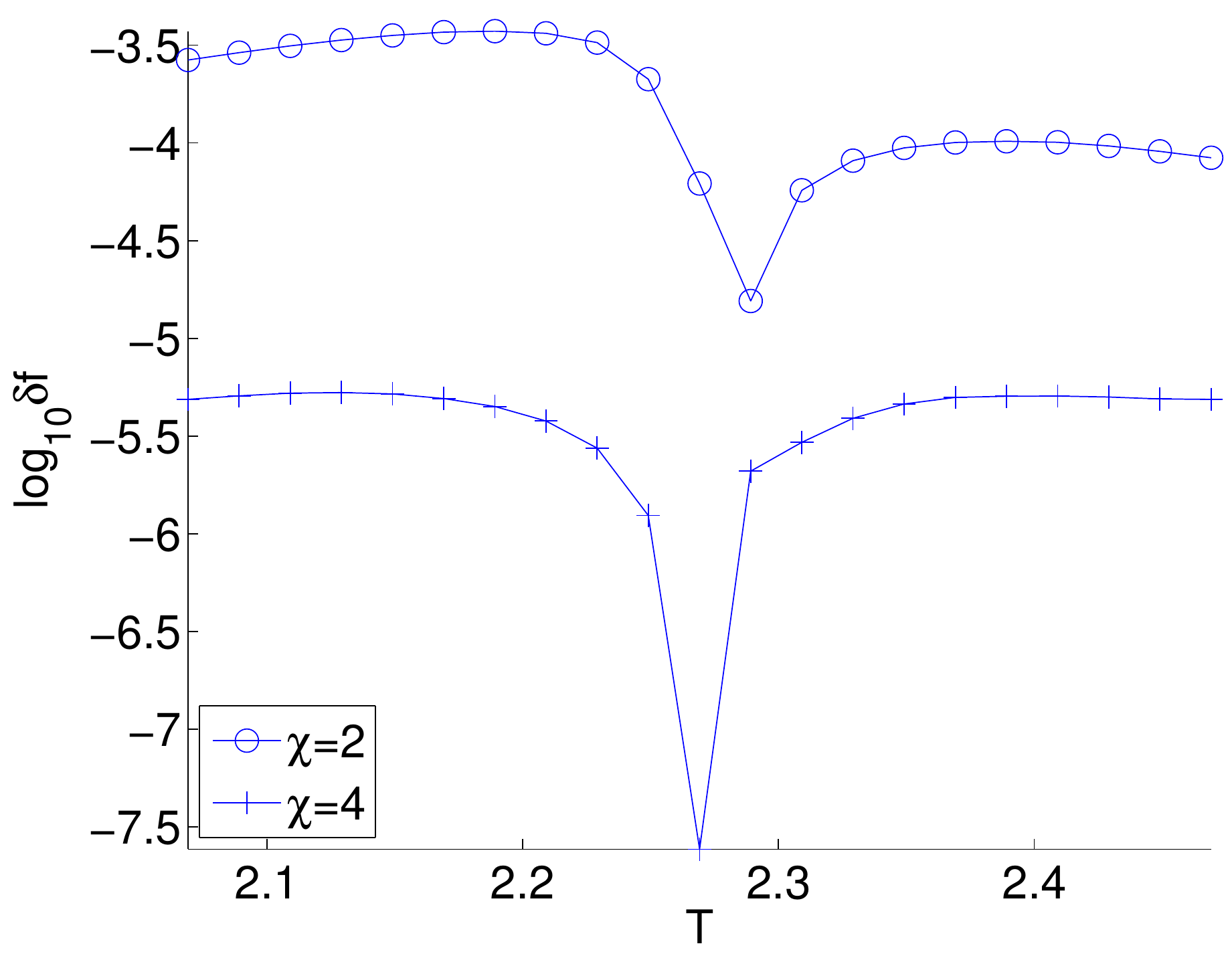} &    \includegraphics[height=2in]{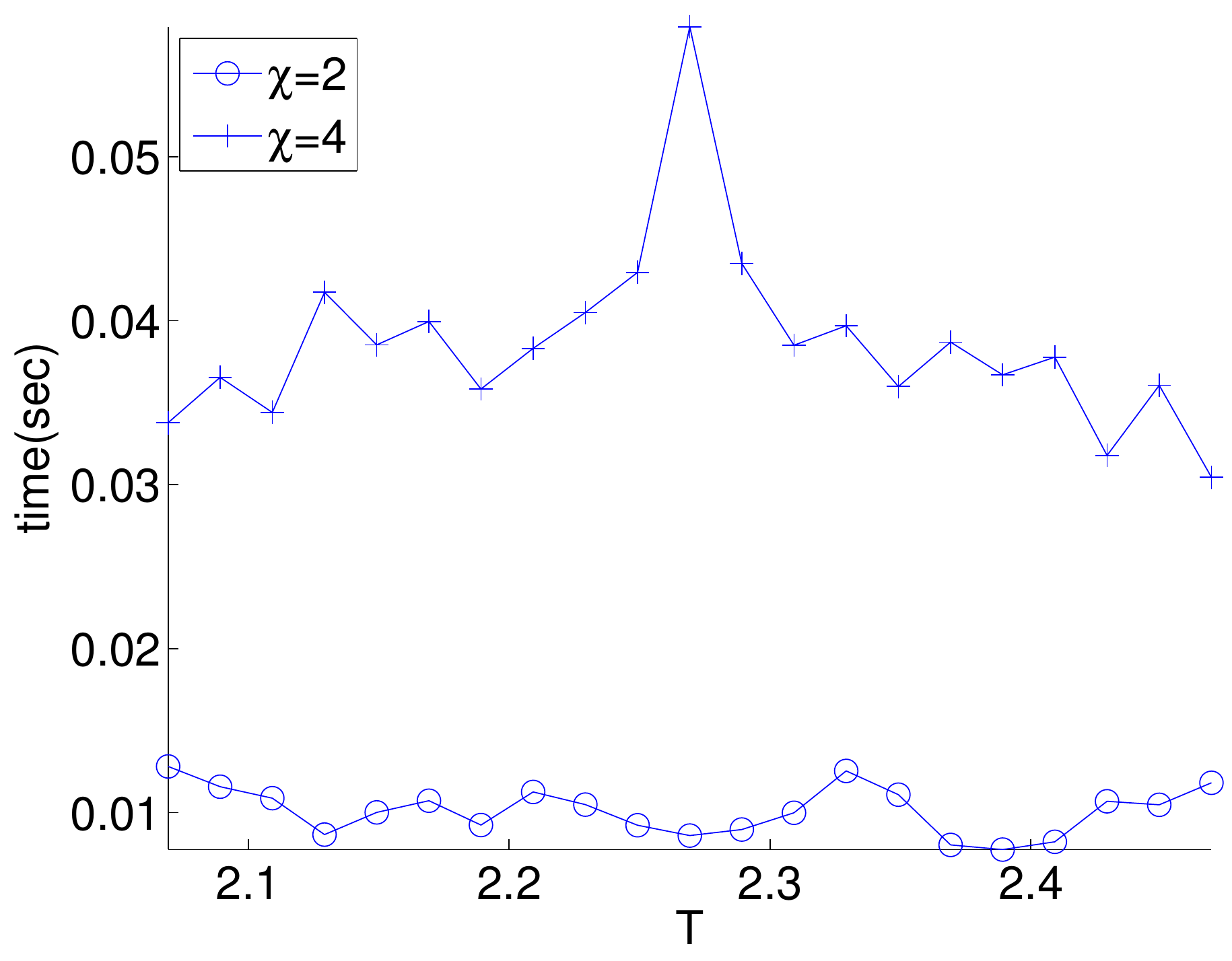}\\
    (a) & (b)
  \end{tabular}
  \caption{Results of free energy calculation. (a) The relative error
    $\delta f_N(\beta)$ of the free energy per site at temperatures
    around $T_c$ for $\chi=2,4$. (b) The running time per iteration
    of TNS for the same $T$ and $\chi$ values.  }
  \label{fig:2Dfree}
\end{figure}

Let us denote by $\tilde{Z}_N(\beta)$ the numerical approximation of
the partition function $Z_N(\beta)$ obtained via TNS. The exact free
energy per site $f_N(\beta)$ and the approximate free energy per site
$\tilde{f}_N(\beta)$ are defined by
\[
f_N(\beta) = \left(-\frac{1}{\beta} \log Z_N(\beta)\right) / N, \quad
\tilde{f}_N(\beta) = \left(-\frac{1}{\beta} \log \tilde{Z}_N(\beta)\right) / N.
\]
For an infinite 2D statistical Ising system, the free energy per site
\[
f(\beta) = \lim_{N\rightarrow\infty} f_N(\beta,N)
\]
can be derived analytically \cite{Huang1987}. Therefore, for
sufficiently large $N$, $f_N(\beta,N)$ is well approximated by
$f(\beta)$. In order to measure the accuracy of TNS for computing the
partition function, we define the relative error
\[
\delta f_N(\beta) 
\equiv \frac{|\tilde{f}_N(\beta)-f(\beta)|}{|f(\beta)|}
\approx \frac{|\tilde{f}_N(\beta)-f_N(\beta)|}{|f_N(\beta)|}.
\]
The critical temperature of the 2D statistical Ising model is
$T_c=1/\ln(1+\sqrt{2})$. For a periodic statistical Ising model on a
$2^{15}\times 2^{15}$ lattice, Figure \ref{fig:2Dfree} plots the
relative error (left) and the running time per iteration (right) for
$\chi=2,4$ at different temperatures near the critical temperature
$T_c$.

From the plots in Figure \ref{fig:2Dfree} one can make the
following observations.
\begin{itemize}
\item First, TNS removes the short-range correlation quite
  effectively. With $\chi=4$, it achieves 5-6 digits of accuracy for
  the relative free energy per site. Even with $\chi=2$, one obtains
  3-4 digits of accuracy.
\item Second, TNS is quite efficiently. For $\chi=4$, each iteration
  of the TNS takes about 0.05 seconds. The running time tends to grow
  a bit when $T$ approaches the critical temperature $T_c$.
\item Most surprisingly, for a fixed $\chi$ value, TNS gives more
  accurate results when the temperature is close to $T_c$. For example
  with $\chi=4$ and at $T=T_c$, the relative error is on the order of
  $10^{-8}$. This is drastically different from most of the TRG-type
  algorithms where the accuracy deteriorates significantly near $T_c$.
\end{itemize}

%-----------------------------------
\subsection{Observables}\label{sec:2Dob}

The TNS algorithm described in Section \ref{sec:2Dpf} for computing
the partition function (and equivalently the free energy) can be
extended to compute observables such as the average magnetization and
the internal energy per site.

The internal energy $U_N(\beta)$ of the whole system and the internal
energy per site $u_N(\beta)$ are defined as
\[
  U_N(\beta) = \p_\beta(-\log Z_N(\beta)) = -\frac{\p_\beta Z_N(\beta)}{Z_N(\beta)},\quad
  u_N(\beta) = \frac{U_N(\beta)}{N}.
\]
A direct calculation shows that 
\[
\p_\beta Z_N(\beta) = \sum_\sigma e^{-\beta H_N(\sigma)}(-H_N(\sigma)) =
\sum_\sigma (\sum_{(ij)} \sigma_i\sigma_j) e^{-\beta H_N(\sigma)} =
N_e \sum_\sigma (\sigma_i\sigma_j) e^{-\beta H_N(\sigma)}
\]
where in the last formula $(i,j)$ can be any edge due to the
translational invariance of the system and $N_e=2N$. This gives the
following formula for the internal energy per site
\begin{equation}\label{eq:u}
u_N(\beta) = \frac{U_N(\beta)}{N} = 
\frac{N_e}{N}\cdot
\frac{\sum_\sigma(\sigma_i\sigma_j) e^{-\beta H_N(\sigma)}}{\sum_\sigma e^{-\beta H_N(\sigma)}}=
2\frac{\sum_\sigma(\sigma_i\sigma_j) e^{-\beta H_N(\sigma)}}{\sum_\sigma e^{-\beta H_N(\sigma)}}.
\end{equation}

To define the average magnetization, one introduces a small external
magnetic field $B$ and defines the partition function of the this
perturbed system
\[
Z_{N,B}(\beta) = \sum_\sigma e^{-\beta H_{N,B}(\sigma)},\quad
H_{N,B}(\sigma)=-\left(\sum_{(ij)}\sigma_i\sigma_j+B\sum_i\sigma_i\right).
\]
The magnetization at a single site $i$ is equal to
\begin{equation}\label{eq:m}
\average{\sigma_i}_{N,B}(\beta) = 
\frac{\sum_\sigma \sigma_i e^{-\beta H_{N,B}(\sigma)}}
     {\sum_\sigma e^{-\beta H_{N,B}(\sigma)}},
\end{equation}
and the average magnetization $m_{N,B}(\beta)$ is equal to the same
quantity since
\[
m_{N,B}(\beta) = \frac{1}{N} \sum_i \average{\sigma_i}_{N,B}(\beta) =
\average{\sigma_i}_{N,B}(\beta),
\]
where in the last formulation $i$ can be any site in the periodic
Ising model due to the translational invariance of the system.

%----
\subsubsection{Algorithm}\label{sec:2Dobal}

The computation of the quantities mentioned above requires the
evaluation of the following sums:
\begin{equation}\label{eq:ob}
  \sum_\sigma(\sigma_i\sigma_j) e^{-\beta H_N(\sigma)},\quad
  \sum_\sigma \sigma_i e^{-\beta H_{N,B}(\sigma)},
\end{equation}
where $i$ is any site in the first formula while $(i,j)$ is any bond
in the second. Both sums can also be represented using tensor
networks using the so-called impurity tensor method. 

%% The main idea is the the tensor networks for these two quantities can
%% be obtained by replacing a single tensor in the tensor network of the
%% partition function.

Recall that the 2D periodic statistical Ising model considered here is
of size $n\times n$ where $n=2^L$. Without loss of generality, one can
assume that the sites $i$ and $j$ in \eqref{eq:ob} are located inside
the $2\times 2$ sub-lattice at the center of the whole computation
domain. Following the same reasoning in Section \ref{sec:intro}, one
can represent $\sum_\sigma(\sigma_i\sigma_j) e^{-\beta H_N(\sigma)}$
and $\sum_\sigma \sigma_i e^{-\beta H_{N,B}(\sigma)}$ as tensor
networks. The only difference between them and the tensor network of
$Z_N(\beta)$ is a single tensor located inside this $2\times 2$
sub-lattice at the center. 

\begin{figure}[h!]
  \includegraphics[scale=0.2]{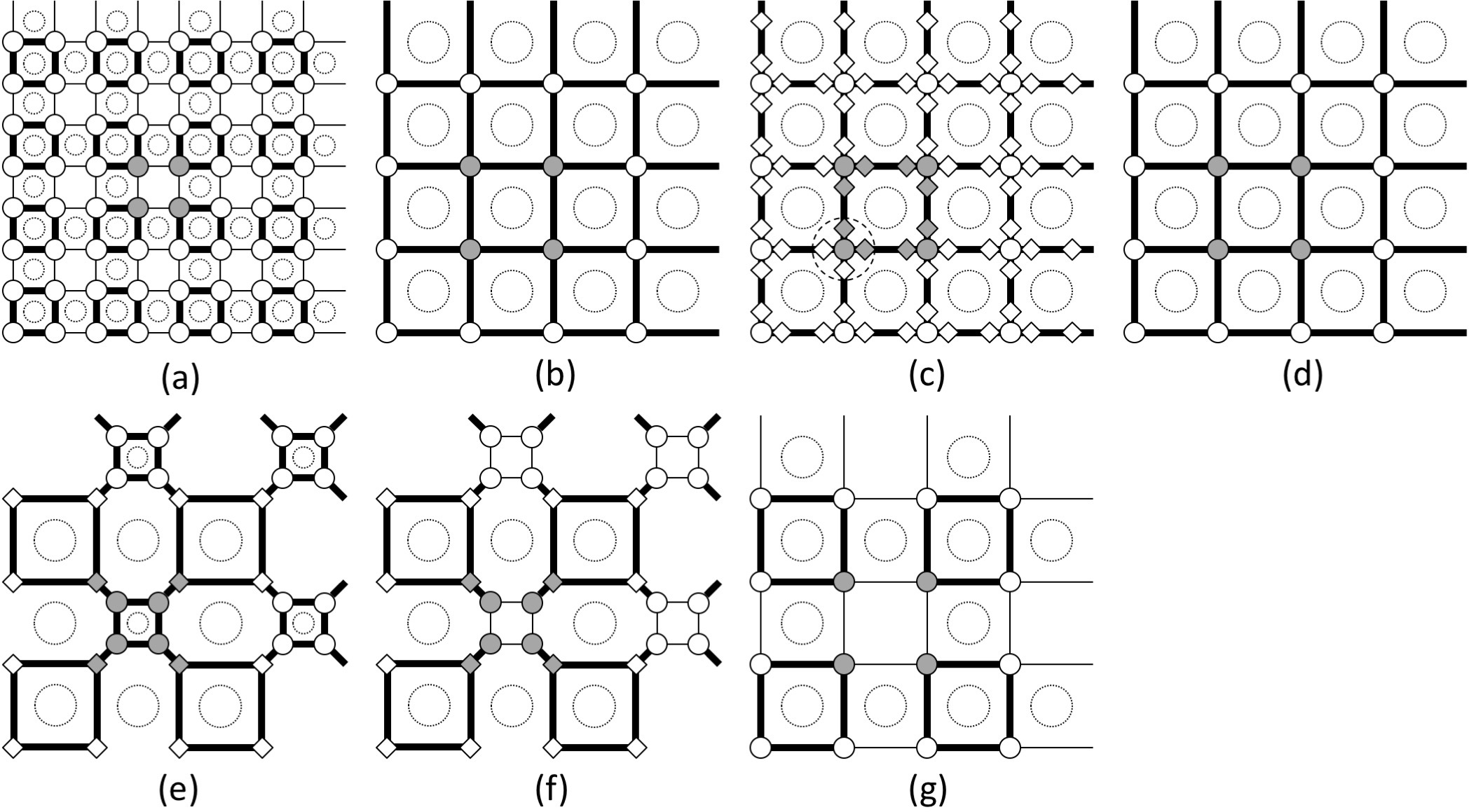}
  \caption{Impurity method for computing the spontaneous magnetization
    and the internal energy per site. The iteration invariance also
    requires that at the beginning of each iteration only the four
    tensors near the center can be different from the ones used in
    $Z_N(\beta)$. These four special tensors are marked in gray. At
    the end of each iteration, one obtains an coarse-grained tensor
    network that also satisfies this condition.  }
  \label{fig:impurity}
\end{figure}

The algorithm for computing these new tensor networks are quite
similar and it becomes particularly simple when the modified TNS
algorithm in Section \ref{sec:2Dpfmod} is used. The whole algorithm is
illustrated in Figure \ref{fig:impurity} and here we highlight the
main differences.
\begin{itemize}
\item In addition to the iteration invariance of the modified
  algorithm in Section \ref{sec:2Dpfmod}, one also requires that only
  the four tensors at the center (marked in gray in Figure
  \ref{fig:impurity}(a)) can be different from the ones used for
  $Z_N(\beta)$.
\item Because the four special tensors are at the center at the tensor
  network at level $\ell$, after contraction there are exactly four
  special tensors at the center of the tensor network at level
  $\ell+1$ (marked in gray in Figure \ref{fig:impurity}(b)). The rest
  are identical to the ones used for $Z_N(\beta)$.
\item From Figure \ref{fig:impurity}(b) to Figure
  \ref{fig:impurity}(c), the $UU'T$-projections at the four
  surrounding edges of the center plaquette are computed from the four
  special corner tensors. The resulting orthogonal $U$ matrices are
  marked in gray as well. The $UU'T$-projection at all other edges are
  inherited from the algorithm for the partition function
  $Z_N(\beta)$. When contracting the tensor at each vertex with its
  four adjacent orthogonal (diamond) matrices (see Figure
  \ref{fig:impurity}(c) to Figure \ref{fig:impurity}(d)), this ensures
  that only the four tensors at the center are different from the ones
  used for $Z_N(\beta)$.
\item In the structure-preserving skeletonization step for the
  $(1,1)_2$ plaquettes (see Figure \ref{fig:impurity}(e) and (f)),
  only the center $(1,1)_2$ plaquette is different from the one
  appeared in $Z_N(\beta)$. Therefore, this is the only one that
  requires an extra structure-preserving skeletonization computation.
\item When contracting the tensors at the corners of the $(1,1)_2$
  plaquettes to get back the 4-tensors in Figure
  \ref{fig:impurity}(g), again only the four tensors at the center
  (marked in gray) are different. This ensures that the tensor network
  at the beginning of the next iteration satisfies the iteration
  invariance mentioned above.
\end{itemize}
At each iteration of the in this impurity method, the algorithm
performs a constant number of extra $UU'T$-projection and one extra
structure-preserving skeletonization for the $(1,1)_2$ plaquette at
the center. When $\chi$ is fixed, all these computation takes a
constant number of steps. As a result,the extra computational cost for
the impurity method is proportional to $O(L)=O(\log N)$ once the
evaluation of $Z_N(\beta)$ is ready.

%----
\subsubsection{Numerical results}

For the internal energy $u_N(\beta)$, we denote by
$\tilde{u}_N(\beta)$ its TNS approximation. When $N$ approaches
infinity, the limit
\[
u(\beta) = \lim_{N\rightarrow \infty} u_N(\beta)
\]
can be derived analytically \cite{Huang1987}. Therefore for $N$
sufficiently large, $u(\beta)$ serves as a good benchmark for measuring
the accuracy of the TNS algorithm.

For the averaged magnetization, let us denote by
$\tilde{m}_{N,B}(\beta)$ the TNS approximation of
$m_{N,B}(\beta)$. For the 2D statistical Ising model, the spontaneous
magnetization $m_+(\beta)$ is defined as
\[
m_+(\beta) = \lim_{B\rightarrow 0^+} \lim_{N\rightarrow \infty} m_{N,B}(\beta)
\]
and this can be written down analytically as well
\cite{Huang1987,Yang1952}. When $B$ is a small positive number, by
setting $N$ to be sufficiently large, one can treat $m_+(\beta)$ as a
good approximation of $m_{N,B}(\beta)$ and use it as a benchmark for
measuring the accuracy of $\tilde{m}_{N,B}(\beta)$.

\begin{figure}[h!]
  \begin{tabular}{cc}
    \includegraphics[height=2in]{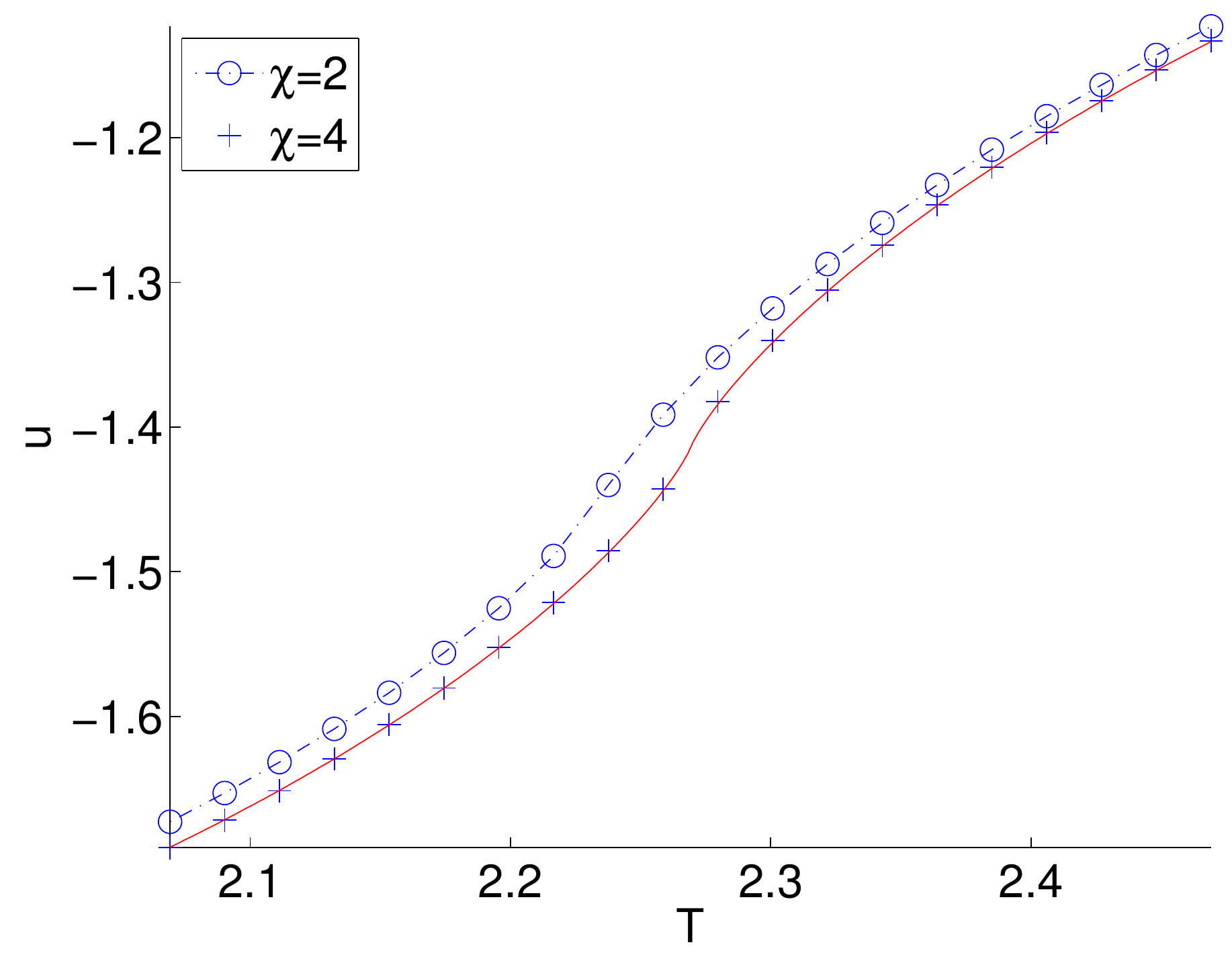}&    \includegraphics[height=2in]{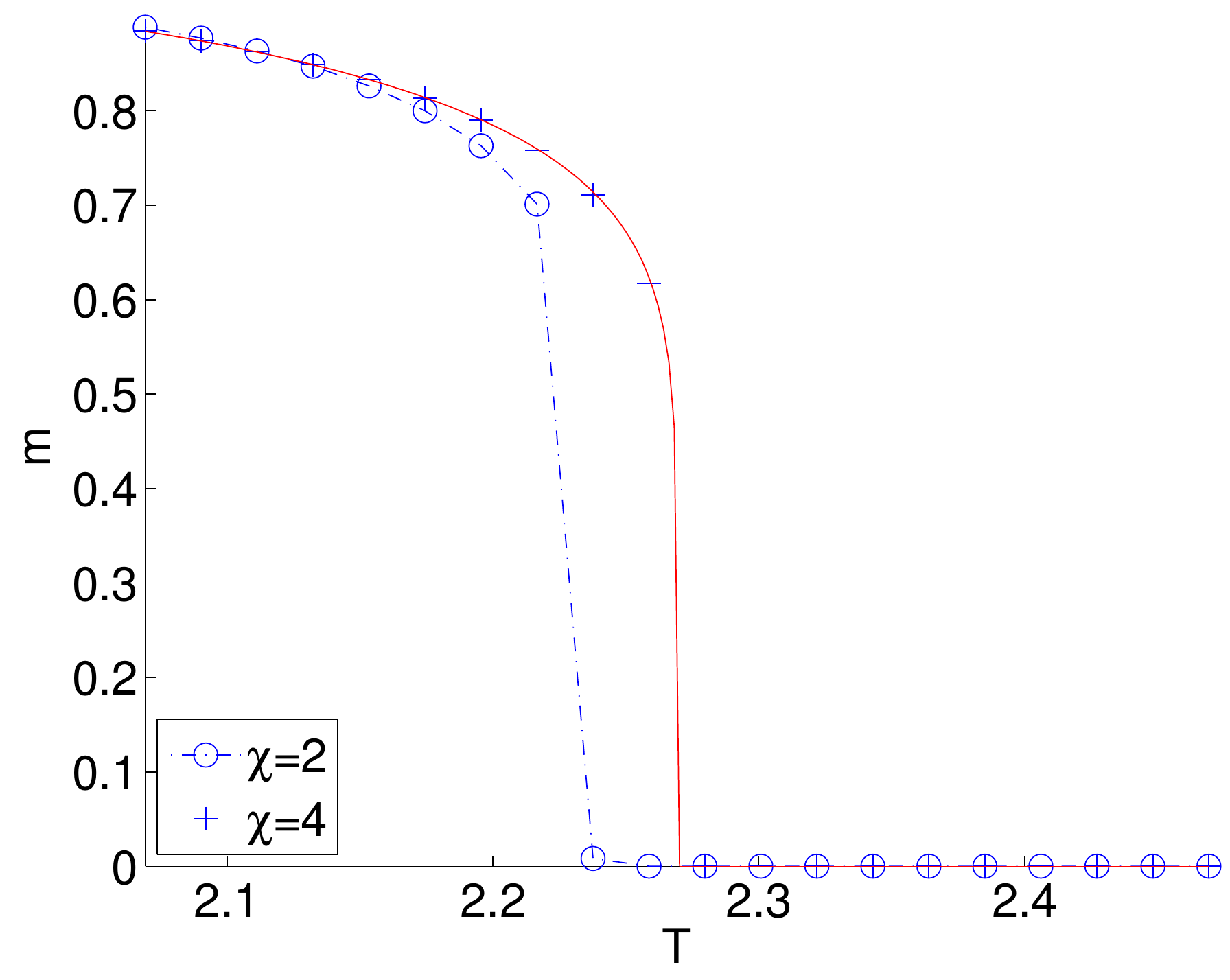} \\
    (a) & (b)
  \end{tabular}
  \caption{Numerical results for computing the observables using the
    impurity method for $\chi=2,4$.  (a) Internal energy. (a) Average
    magnetization. }
  \label{fig:magnetic}
\end{figure}

Figure \ref{fig:magnetic}(a) shows the computed internal energy per
site $\tilde{u}_N(\beta)$ along with $u(\beta)$ for $\chi=2,4$. On the
right, Figure \ref{fig:magnetic}(b) gives the computed average
magnetization $\tilde{m}_{N,B}(\beta)$ along with the spontaneous
magnetization $m_+(\beta)$. Though the computation with $\chi=2$ has a
significant error, it does exhibit the phase-transition clearly. Once
$\chi$ is increased to 4, the numerical results and the exact curves
match very well.

%--------------
\subsection{Extension to disordered systems}

The tensor network skeletonization algorithm can also be extended
easily to disordered systems and we briefly sketch how this can be
done. For example, consider the 2D Edwards-Anderson spin-glass model
(see \cite{Nishimori2001} for example) where the spins are arranged
geometrically the same fashion as the classical Ising model but each
edge $(i,j)$ is associated with a parameter $J_{ij}$. For a fixed
realization of $J\equiv \{J_{ij}\}$, the partition function is given
by
\[
Z_{N,J}(\beta) = \sum_\sigma e^{-\beta H_{N,J}(\sigma)},\quad
H_{N,J}(\sigma) = -\sum_{(ij)}J_{ij}\sigma_i\sigma_j.
\]
At a fixed realization of $J_{ij}$, the order parameter of the model
is defined as
\[
q_{N,J}(\beta)=\frac{1}{N}\sum_{i=1}^N \average{\sigma_i}_{N,J}^2(\beta),\quad
\average{\sigma_i}_{N,J}(\beta) = \frac{\sum_\sigma \sigma_i e^{-\beta
    H_{N,J}(\sigma)}}{\sum_\sigma e^{-\beta H_{N,J}(\sigma)}}.
\]

The computation of the order parameter $q_N(\beta)$ first requires the
evaluation of $Z_{N,J}(\beta)$. Similar to the standard Ising model,
this can be represented with a tensor network. The TNS algorithm
remains the same, except that the computation at each plaquette has to
be performed separately since the system is not translationally
invariant anymore. For any fixed bond dimension $\chi$, the
computational complexity of TNS scales like $O(N)$, where $N$ is the
number of spins.

It also requires the evaluation of $\sum_\sigma \sigma_i e^{-\beta
  H_J(\sigma)}$ for each $i$. The discussion in Section \ref{sec:2Dob}
shows that for each $i$ one needs to perform $O(\log N)$ extra
structure-preserving skeletonizations, since most of the computation
of $Z_{N,J}(\beta)$ can be reused. Therefore, the computation of
$\average{\sigma_i}_{N,J}^2(\beta)$ for all spins $i$ takes $O(N\log
N)$ steps. Putting this and the cost of evaluating $Z_{N,J}(\beta)$
together shows that the computation of the order parameter
$q_{N,J}(\beta)$ can be carried out in $O(N\log N)$ steps.

%-----------------------------------
\section{TNS for 3D statistical Ising model}\label{sec:3D}

In this section, we describe how to extend the tensor network
skeletonization algorithm to the 3D statistical Ising model. One key
feature that has not been emphasized is that TNS preserves the
Cartesian structure of the problem. This allow for a natural
generalization to 3D systems. Let us consider a 3D periodic
statistical Ising model on an $n\times n\times n$ Cartesian
grid. $N=n^3$ is the number of total spins and we assume without loss
of generality $n=2^L$ for an integer $L$.

%--------------
\subsection{Partition function}\label{sec:3Dpf}
The partition function $Z_N(\beta)$ can be represented with a tensor
network $(V^0, E^0, \{T^i\}_{i\in V^0})$ where $V^0$ is the set of
vertices of the Cartesian grid, the edge set $E^0$ contains the edges
between two adjacent sites in the $x$, $y$ and $z$ directions, and
$T^i$ is a 6-tensor at site $i$. This gives rise to an $n\times
n\times n$ array of small cubes, each with its 8 vertices in $V_0$. If
a cube has vertex $i=(i_1,i_2,i_3)$ at its lower-left-front corner,
then we shall index this cube with $i$ as well. We refer to the cubes
with index equal to $(0,0,0)$ modulus 2 as $(0,0,0)_2$ cubes and those
with index equal to $(1,1,1)$ modulus 2 as $(1,1,1)_2$ cubes. As with
the 2D case, we let $\chi$ be a predetermined upper bound for the bond
dimension and without loss of generality one can assume that $\chi_e
\approx \chi$ for each $e\in E^0$.

%------
\subsubsection{Algorithm}\label{sec:3Dpfal}

\begin{figure}[h!]
  \includegraphics[scale=0.2]{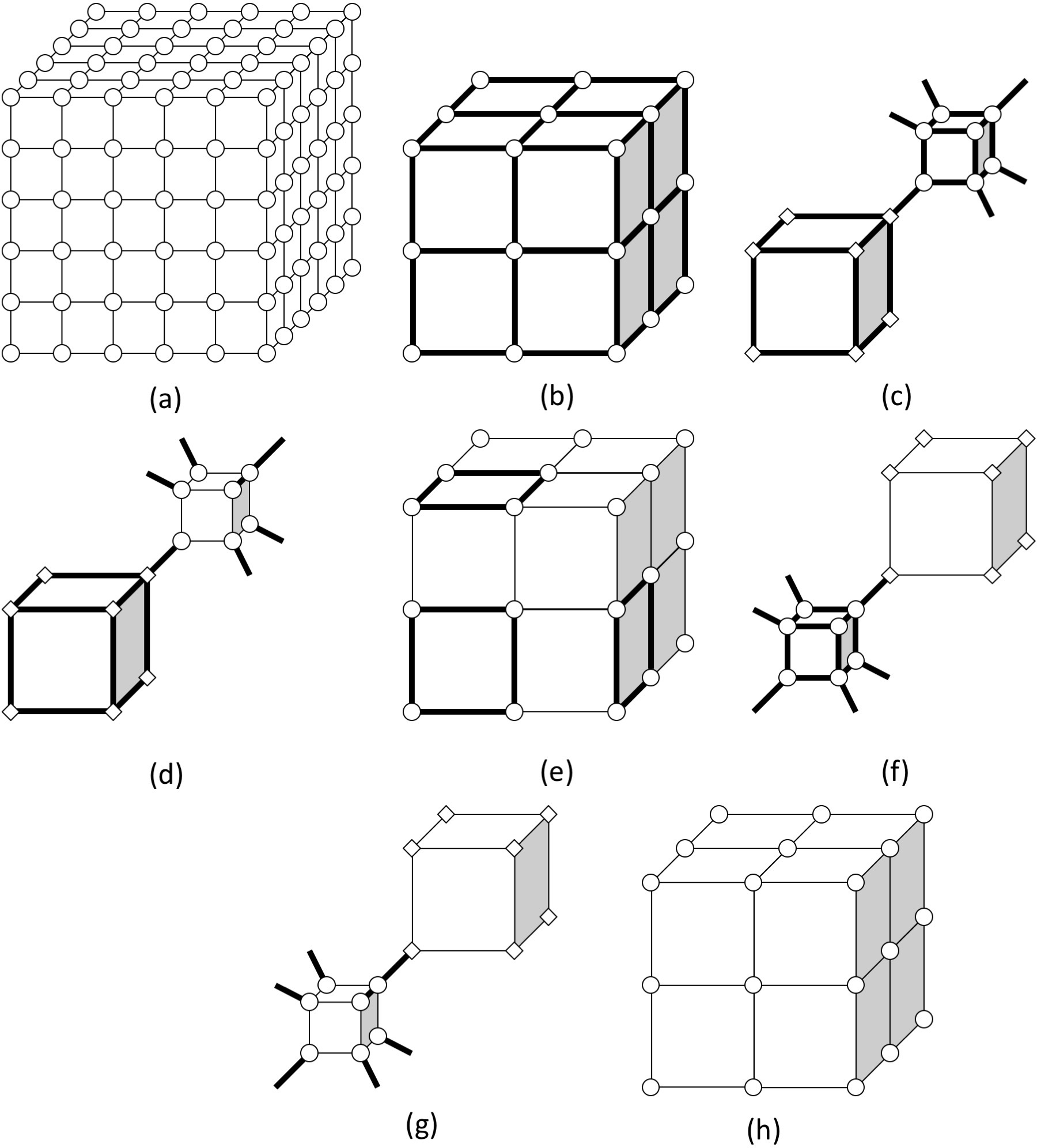}
  \caption{A single iteration of the tensor network skeletonization
    algorithm. The starting point is a tensor network with bond
    dimension $\chi$ and short-range correlation removed in $(0,0,0)_2$
    and $(1,1,1)_2$ cubes. The final point is a coarse-grained
    tensor-network with 1/8 vertices (tensors). This coarse-grained tensor
    also has bond dimensions equal to $\chi$ and has short-range
    correlation removed for the (larger) $(0,0,0)_2$ and $(1,1,1)_2$
    cubes.
  }
  \label{fig:3D}
\end{figure}

The TNS algorithm consists of a sequence of coarse-graining
iterations. At the beginning of each iteration (except the $0$-th
iteration), we require the following iteration invariance to hold:
\begin{itemize}
\item for each of the $(0,0,0)_2$ and $(1,1,1)_2$ cubes, the
  short-range correlation has already been eliminated.
\end{itemize}

At the beginning of the $\ell$-th iteration, one holds a tensor
network $(V^\ell,E^\ell,\{T^i\}_{i\in V_\ell})$ at level $\ell$ with
$(n/2^\ell)\times (n/2^\ell) \times (n/2^\ell)$ vertices. The
$\ell$-th iteration consists of the following steps.
\begin{enumerate}
\item 
  Contract the tensors at the eight vertices of each $(0,0,0)_2$ cube
  into a single tensor (see Figure \ref{fig:3D}(b)). The $(1,1,1)_2$
  cubes are stretched and this results a new tensor network that
  contains $1/8$ of the vertices. The tensors at the new vertices are
  identical and the bond dimension of the new edges are equal to
  $\chi^4$ (shown with bold lines in the figure).  Similar to the 2D
  case, the short-range correlations at level $\ell$ does not survive
  to level $\ell+1$ due to the iteration invariance. However, there
  are short-range correlations for the cubes at level $\ell+1$. The
  key task is to remove some of these short-range correlations and
  reduce the bond dimension back to $\chi$.
\item At each vertex $i$ in Figure \ref{fig:3D}(b), denote the tensor
  by $T^i_{abcdef}$ where $a$, $b$, $c$, $d$, $e$, and $f$ are the
  left, right, bottom, top, front, and back edges, respectively. By
  invoking three $UU'T$-projection step (one for each of the left,
  bottom, and front edges), one effectively inserts two orthogonal
  (diamond) matrices in each of these edges. At each vertex $i$,
  further merge the tensor $T^i_{abcdef}$ with the six adjacent
  orthogonal (diamond) matrices. This step does not change the topology
  of the tensor network but the $T^i$ tensor has been modified.
\item For each $(1,1,1)_2$ cube in Figure \ref{fig:3D}(b), apply the
  $UR$-projection to the 6-tensor at each of its corners. Here the
  three edges adjacent to the cube are grouped together. Notice that
  the round $R$ tensors are placed close to the $(1,1,1)_2$ cube. This
  projection step only keeps the top $\chi^3$ singular values, i.e.,
  the bond dimension of the diagonal edges are equal to $\chi^3$. The
  resulting graph is given in Figure \ref{fig:3D}(c).
\item In this key step, apply structure-preserving skeletonization to
  each of the $(1,1,1)_2$ cubes. The details of this procedure will be
  provided in Section \ref{sec:3Dpfsps}. The resulting $(1,1,1)_2$
  plaquette has short-range correlation removed and the bond
  dimensions of its 12 surrounding edges are reduced from $\chi^4$ to
  $\chi$ (see Figure \ref{fig:3D}(d)). One then merges back the
  $UR$-projections at its eight corners. The resulting tensor network
  in Figure \ref{fig:3D}(e) is similar to the one in Figure
  \ref{fig:3D}(b) but the short-range correlations in the $(1,1,1)_2$
  plaquettes are now removed.
\item Now repeat the previous two steps to the $(0,0,0)_2$
  plaquettes. This is illustrated in Figure \ref{fig:3D}(f), (g) and
  (h). The resulting tensor network has short-range correlation
  removed in both $(0,0,0)_2$ and $(1,1,1)_2$ plaquettes and the bond
  dimension of the edges is reduced back to $\chi$ from $\chi^4$.
\end{enumerate}
This finishes the $\ell$-th iteration. At this point, one obtains a
new tensor network denoted by $(V^{\ell+1},E^{\ell+1},\{T^i\}_{i\in
  V_{\ell+1}})$ that is a self-similar and coarse-grained version of
$(V^\ell,E^\ell,\{T^i\}_{i\in V_\ell})$. This network satisfies the
iteration invariance and can serve as the starting point of the next.
iteration of the algorithm.

The last tensor network $(V^L,E^L,\{T^i\}_{i\in V_L})$ contains only a
single $6$-tensor with the left and right edges identified and
similarly for the bottom/top edges and front/back edges. Contracting
this final tensor gives an approximation for the partition
function. Similar to the 2D case, one can also introduce a modified
version of this algorithm by removing short-range correlation for the
$(1,1,1)_2$ cubes.
%------
\subsubsection{Structure-Preserving skeletonization}\label{sec:3Dpfsps}

\begin{figure}[h!]
  \includegraphics[scale=0.2]{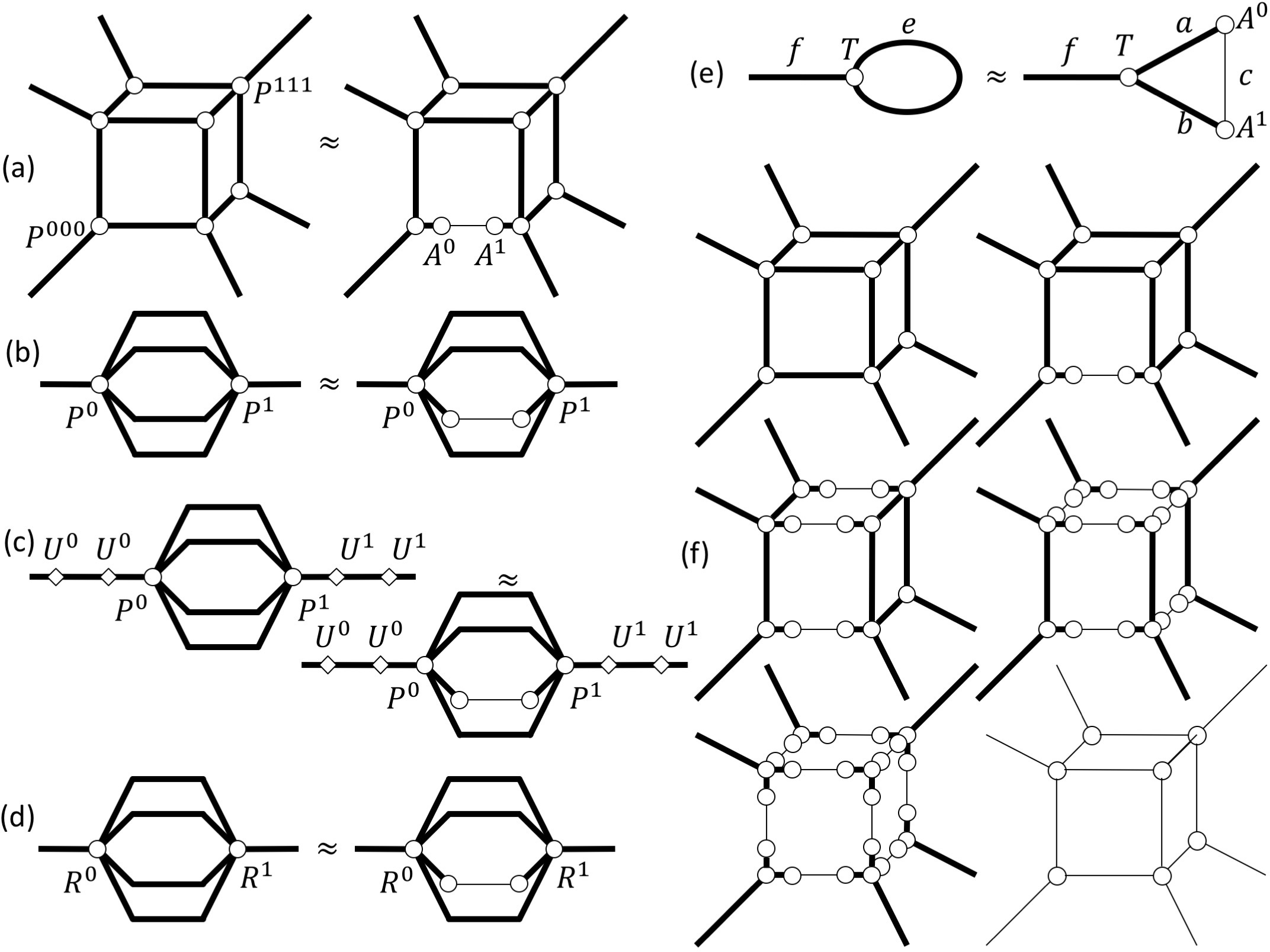}
  \caption{The structure-preserving skeletonization removes the
    short-range correlation within a $(1,1,1)_2$ or $(0,0,0)_2$ cube.}
  \label{fig:3Dskl}
\end{figure}

The structure-preserving skeletonization procedure for the 3D cubes is
similar to the one introduced for 2D plaquette in Section
\ref{sec:2Dpfsps}. This procedure is illustrated in Figure
\ref{fig:3Dskl} with the eight corner 4-tensors denoted by $P^{000}$,
$P^{100}$, $P^{010}$, $P^{110}$, $P^{001}$, $P^{101}$, $P^{011}$, and
$P^{111}$.

Instead of replacing the eight corner 4-tensors of the gray cube
simultaneously, this procedure considers the 12 interior edges one by
one and inserts within each edge two tensors of size $\chi^4 \times
\chi$.
\begin{enumerate}
\item Starting from the bottom front edge, the procedure seeks two
  2-tensors $A^0$ and $A^1$ of size $\chi^4 \times \chi$ with the
  condition that the 8-tensor of the new $(1,1,1)_2$ cube after the
  insertion approximates the original 8-tensors (see Figure
  \ref{fig:3Dskl}(a)).
\item Merge the two left tensors $P^{000}$, $P^{001}$, $P^{010}$, and
  $P^{011}$ into a 5-tensor $P^0$ and merge the four right tensors
  into a 5-tensor $P^1$. After that, the condition is equivalent to
  the one given in Figure \ref{fig:3Dskl}(b) with the two boundary
  edges have bond dimension equal to $(\chi^3)^4 = \chi^{12}$.
\item Since the bond dimensions of the two edges between $P^0$ and
  $P^1$ are to be reduced to $\chi$, this implies that the bond
  dimensions of the two boundary edges can be reduced to $\chi^4$
  instead of $\chi^{12}$.  As a result, one can perform the
  $UU'T$-projection to both $P^0$ and $P^1$. This gives rise the
  condition in Figure \ref{fig:3Dskl}(c).
\item Remove the two tensors $U^0$ and $U^1$ at the two endpoints,
  contract $U^0$ with $P^0$ to get a 3-tensor $R^0$, and contract
  $U^1$ with $P^1$ to get $R^1$. The approximation condition can now
  be written in terms of $R^0$ and $R^1$ as in Figure
  \ref{fig:3Dskl}(d).
\item Finally, contracting the three other edge between $R^0$ and
  $R^1$ results a new 3-tensor $T$. The approximation condition now
  takes the form given in Figure \ref{fig:3Dskl}(e). This is now
  exactly the setting of the skeletonization procedure and can be
  solved using the alternating least square algorithm proposed 
  in Section \ref{sec:basicsps}.
\end{enumerate}
At this point, two tensors $A^0$ and $A^1$ are successfully inserted
into the bottom front edge. One can repeat this process also for the
other three edges in the $x$ direction. Once this is done, we repeat
this for the edges in the $y$ direction and then for the edges in the
$z$ direction. At this point, there are in total 24 orthogonal tensors
inserted in the 12 surrounding edges of the cube. Finally, merging
each of the corner tensors with its three adjacent tensors gives the
desired approximation (see Figure \ref{fig:3Dskl}(f) for the whole
process).

%------
\subsubsection{Numerical results}

\begin{figure}[h!]
  \includegraphics[height=2in]{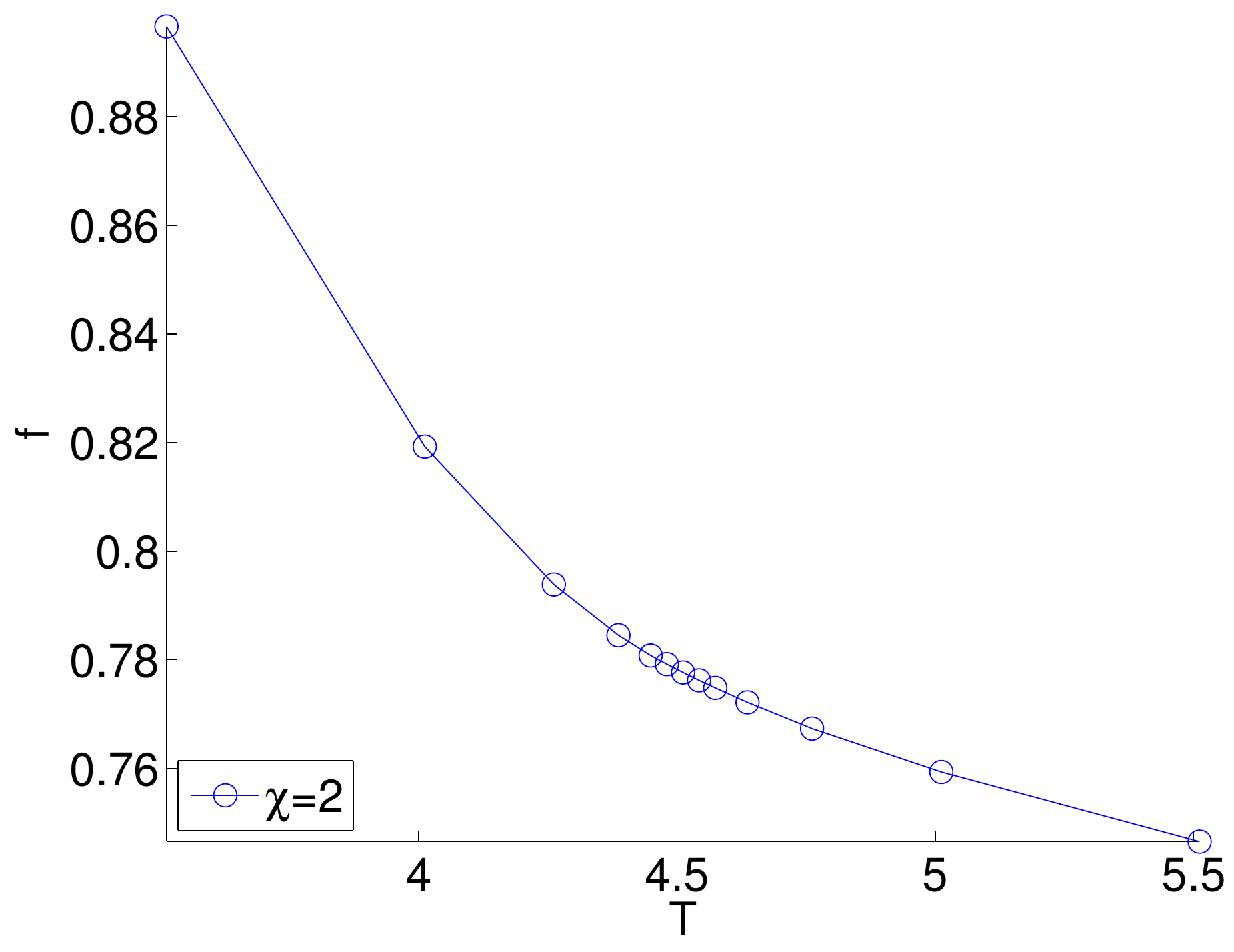}
  \caption{The free energy per site obtained through TNS at
    temperatures around $T_c$ for $\chi=2$ for the 3D periodic
    statistical Ising model.}
  \label{fig:3Dfree}
\end{figure}

The critical temperature of the 3D statistical Ising model is
$T_c\approx 4.5115$ but the free energy per site is not known
explicitly. For a 3D periodic Ising model on a $2^{6}\times 2^{6}$
lattice, Figure \ref{fig:3Dfree} shows the free energy per site
obtained through TNS for $\chi=2$ at different temperatures near the
$T_c$. The obtained values of the free energy is close to the results
obtained from other calculations using HOTRG or Monte Carlo
calculations.

%% In the following figure, we plot the magnetic moments and internal
%% energy.
%% \begin{figure}[h!]
%%   \includegraphics[scale=0.2]{fig/tm_Slide13.jpg}
%%   \caption{3D magnetic/internal}
%%   \label{fig:3Dmagnetic}
%% \end{figure}

%--------------
\subsection{Extensions}

Similar to the 2D case, the 3D algorithm can be used to compute the
average magnetization and the internal energy per site. When
representing these quantities through tensor networks, one finds that
only the 8 tensors at the center of the computational domain are
different from the ones used in the partition function
calculation. Therefore, the impurity tensor method can be applied as
expected and the extra computational cost grows like $O(\log N)$ for
any fixed $\chi$.

For disordered systems, the same discussion for the 2D systems
applies. For example, for computing the order parameter of the 3D
Edwards-Anderson model, one only needs to perform one impurity tensor
computation for each site and thus the overall complexity grows like
$O(N\log N)$ for any fixed $\chi$.

%-----------------------------------
\section{Ground state for quantum Ising models}\label{sec:ground}

\begin{figure}[h!]
  \includegraphics[scale=0.2]{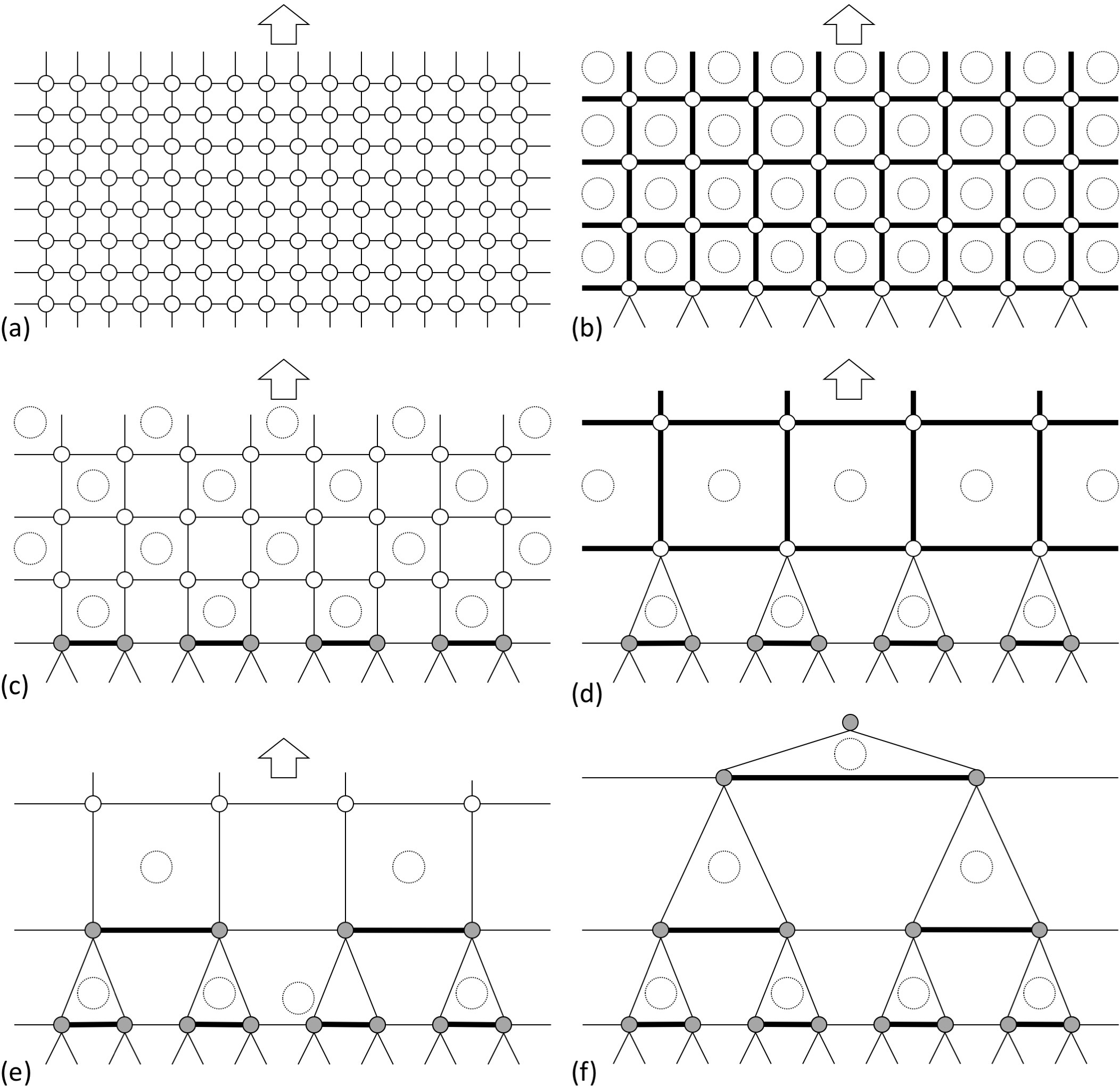}
  \caption{When applied to a Euclidean path integral formulation, TNS
    yields a new representation of the ground state of 1D quantum
    Ising model.}
  \label{fig:1DGS}
\end{figure}

In this section, we briefly touch on how to use TNS to efficiently
represent the ground state of quantum many body system with periodic
boundary condition. Consider for example a 1D periodic quantum Ising
model. One can represent the ground state up to a constant factor
using the Euclidean path integrals \cite{Evenbly2015C}. After some
preliminary tensor manipulations, this turns into a tensor network
that is periodic in the spatial dimension and semi-infinite in the
imaginary time dimension (see Figure \ref{fig:1DGS}(a)).

In the same fashion that the tensor network renormalization (TNR)
gives rise to the multi-scale entanglement renormalization ansatz
(MERA) \cite{Evenbly2015B} for the ground state, TNS generates a new
representation of the ground state as well. Illustrated in Figure
\ref{fig:1DGS}, this process consists of the following steps.
\begin{enumerate}
\item First, contract each group of $2\times 2$ tensors (see Figure
  \ref{fig:1DGS}(b)). The new edges marked with bold lines have bond
  dimension equal to $\chi^2$.
\item Perform the structure-preserving skeletonization to all
  $(0,0)_2$ and $(1,1)_2$ plaquettes to remove the short-range
  correlations and reduce the bond dimension back to $\chi$.  Notice
  that, after the structure-preserving skeletonization, the resulting
  tensors at the bottom level are different from the ones above due to
  their adjacency to the boundary. These special bottom level tensors
  are marked in gray (see Figure \ref{fig:1DGS}(c)).
\item Repeat this process to the remaining tensors above the bottom
  level. Contracting each group of $2\times 2$ tensors results a
  tensor network illustrated in Figure \ref{fig:1DGS}(d) and Figure
  \ref{fig:1DGS}(e). 
\item One can repeat this process until reaching a half-infinite
  string of identical matrices. By extracting its top eigenvector, one
  can reduce this (up to a constant factor) to a 1-tensor at the top
  (see Figure \ref{fig:1DGS}(f)).
\end{enumerate}
The final product is a hierarchical structure shown in Figure
\ref{fig:1DGS}(f)). Though somewhat different from MERA, this new
structure also has the capability of representing strongly entangled
1D quantum systems.

\begin{figure}[h!]
  \includegraphics[scale=0.2]{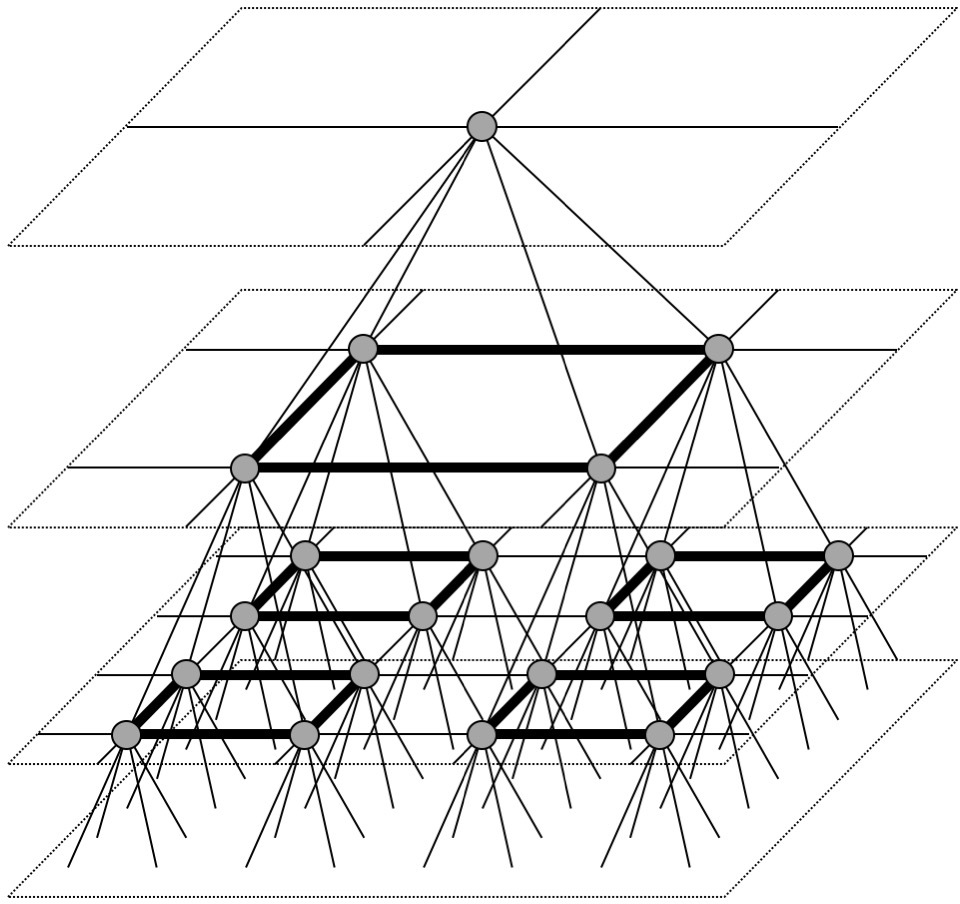}
  \caption{When applied to a Euclidean path integral formulation, TNS
    yields a new representation of the ground state of 2D quantum
    Ising model.}
  \label{fig:2DGS}
\end{figure}

For 2D periodic quantum Ising model, the ground state can be
represented via Euclidean path integral with a 3D tensor network which
is periodic in the $x$ and $y$ directions but semi-infinite in the
imaginary time direction. The above algorithm (with necessary
modifications for the 3D TNS) can be applied to this tensor network
and the result is a hierarchical structure (shown in Figure
\ref{fig:2DGS}) that is capable of representing the ground state of
strongly entangled 2D quantum systems effectively.

%-----------------------------------
\section{Conclusion} \label{sec:conc}

This paper introduced the tensor network skeletonization (TNS) as a
new coarse-graining process for the numerical computation of tensor
networks. At the heart of TNS is a new structure-preserving
skeletonization procedure that removes short-range correlation
effectively.

As to future work, an immediate task is to investigate other
algorithms for the structure-preserving skeletonization problem
\eqref{eq:sps} and \eqref{eq:spsmat}. The alternating least square
algorithm adopted here works quite well in practice. However, it would
be interesting to understand why and also to consider other
alternatives without using the somewhat artificial regularization
parameter.

Most TNS algorithms introduced here are presented in their simplest
forms in order to illustrate the main ideas. This means that they are
not necessarily the most efficient implementations in practice. For
example in the TNS algorithm for partition functions, one performs the
contractions over all directions first and then applies the
$UU'T$-projections to these directions. However in practice, it is
much more efficient to iterate over the directions and, for each
direction, apply a $UU'T$-projection right after the contraction of
this direction.

We also plan to improve on the current implementations for the
disordered systems and the ground state computations as well.

\bibliographystyle{abbrv} \bibliography{ref}

\end{document}